\documentclass{icmart}
\usepackage{amscd}

\contact[braval@math.brown.edu]{Department of Mathematics, Brown University, Providence RI }

\newtheorem{theorem}{Theorem}[section]

\newtheorem{lemma}[theorem]{Lemma}

\theoremstyle{definition}

\newtheorem{remark}[theorem]{Remark}

\theoremstyle{plain}
\newtheorem{Thm}[subsection]{Theorem}
\newtheorem{Cor}[subsection]{Corollary}

\newtheorem{Prop}[subsection]{Proposition}
\newtheorem{Conj}[subsection]{Conjecture}

\theoremstyle{definition}
\newtheorem{Def}[subsection]{Definition}

\theoremstyle{remark}

\numberwithin{equation}{section}
\renewcommand{\rm}{\normalshape}

\newif\ifShowLabels
\ShowLabelstrue
\newdimen\theight
\def\TeXref#1{%
    \leavevmode\vadjust{\setbox0=\hbox{{\tt
        \quad\quad  {\small \rm #1}}}%
    \theight=\ht0
    \advance\theight by \lineskip
    \kern -\theight \vbox to
    \theight{\rightline{\rlap{\box0}}%
    \vss}%
    }}%

\ShowLabelsfalse% comment this out if labels should be printed

\newcommand{\ssec}[2]{\subsection{#2}\label{SS:#1}%
    \ifShowLabels \TeXref{{SS:#1}} \fi}

\newcommand{\refss}[1]{Section ~\ref{SS:#1}}

\newcommand{\refe}[1]{\eqref{E:#1}}

\newenvironment{thm}[1]%
    { \begin{Thm} \label{T:#1}  \ifShowLabels \TeXref{T:#1} \fi }%
    { \end{Thm} }

\renewcommand{\th}[1]{\begin{thm}{#1} \sl }
\renewcommand{\eth}{\end{thm} }

\newcommand{\lem}[1]{\begin{lemma}{#1} \sl}
\newcommand{\elem}{\end{lemma}}

\newenvironment{propos}[1]%
    { \begin{Prop} \label{P:#1}  \ifShowLabels \TeXref{P:#1} \fi }%
    { \end{Prop} }
\newcommand{\prop}[1]{\begin{propos}{#1}\sl }
\newcommand{\eprop}{\end{propos}}

\newenvironment{corol}[1]%
    { \begin{Cor} \label{C:#1}  \ifShowLabels \TeXref{C:#1} \fi }%
    { \end{Cor} }
\newcommand{\cor}[1]{\begin{corol}{#1} \sl }
\newcommand{\ecor}{\end{corol}}

\newenvironment{conjec}[1]%
    { \begin{Conj} \label{Co:#1}  \ifShowLabels \TeXref{Co:#1} \fi }%
    { \end{Conj} }
\newcommand{\conj}[1]{\begin{conjec}{#1} \sl }
\newcommand{\econj}{\end{conjec}}

\newenvironment{defeni}[1]%
    { \begin{Def} \label{D:#1}  \ifShowLabels \TeXref{D:#1} \fi }%
    { \end{Def} }
\newcommand{\defe}[1]{\begin{defeni}{#1} \sl }
\newcommand{\edefe}{\end{defeni}}

\newcommand{\rem}[1]{\begin{remark}{#1}}
\newcommand{\erem}{\end{remark}}

\newcommand{\eq}[1]%
    { \ifShowLabels \TeXref{E:#1} \fi
       \begin{equation} \label{E:#1} }
\newcommand{\eeq}{ \end{equation} }

\newcommand{\prf}{ \begin{proof} }
\newcommand{\epr}{ \end{proof} }
\newcommand\gam{\gamma}     
     
\newcommand\eps{\varepsilon}

\newcommand\kap{\kappa}
\newcommand\lam{\lambda}        \newcommand\Lam{\Lambda}

\newcommand\ome{\omega}     

\newcommand\calA{{\mathcal{A}}}

\newcommand\calF{{\mathcal{F}}}

\newcommand\calH{{\mathcal{H}}}

\newcommand\calK{{\mathcal{K}}}
\newcommand\calL{{\mathcal{L}}}
\newcommand\calM{{\mathcal{M}}}

\newcommand\calO{{\mathcal{O}}}

\newcommand\calU{{\mathcal{U}}}
\newcommand\calV{{\mathcal{V}}}

\newcommand\calZ{{\mathcal{Z}}}

        \newcommand\bfI{{\mathbf I}}

        \newcommand\bfS{{\mathbf S}}

        \newcommand\bfX{{\mathbf X}}

\newcommand\OO{\mathbb{O}}
\newcommand\PP{\mathbb{P}}
\renewcommand\AA{\mathbb{A}}
\renewcommand\SS{\mathbb{S}}
\newcommand\DD{\mathbb{D}}
\newcommand\FF{\mathbb{F}}
\newcommand\GG{\mathbb{G}}

\newcommand\ZZ{\mathbb{Z}}

\newcommand\CC{\mathbb{C}}

 \newcommand\grg{{\mathfrak{g}}}
 \newcommand\grh{{\mathfrak{h}}}

 \newcommand\grl{{\mathfrak{l}}}
 \newcommand\grm{{\mathfrak{m}}}
 \newcommand\grn{{\mathfrak{n}}}

 \newcommand\grq{{\mathfrak{q}}}

 \newcommand\grt{{\mathfrak{t}}}
 \newcommand\gru{{\mathfrak{u}}}

\newcommand\sdp{\times \hskip -0.3em {\raise 0.3ex
\hbox{$\scriptscriptstyle |$}}} % semidirect product

\newcommand\Aut{\operatorname{Aut}}

\newcommand\Conv{\operatorname{Conv}}

\newcommand\End{\operatorname{End\,}}
\newcommand\Ext{\operatorname{Ext}}

\newcommand\GL{\operatorname{GL}}
\newcommand\Gr{\operatorname{Gr}}
\newcommand\gl{{\bf gl}}

\newcommand\Hom{\operatorname {Hom}}

\newcommand\id{\operatorname{id}}
\newcommand\Id{\operatorname{Id}}

\newcommand\Ker{\operatorname{Ker}}

\newcommand\Perv{\operatorname{Perv}}

\newcommand\Spec{\operatorname{Spec}}

\newcommand\Sym{\operatorname{Sym}}

\newcommand\oS{{\overline{S}}}

\newcommand\oX{{\overline{X}}}

\newcommand\uj{{\underline{j}}}

\newcommand\x{\times}
\newcommand\ten{\otimes}

\newcommand{\ra}{\rangle}
\newcommand{\la}{\langle}

\newcommand\bt{\text{b}}

\newcommand\Bun{\operatorname{Bun}}
\newcommand\oBun{\overline{\Bun}}

\newcommand\IC{\operatorname{IC}}
\newcommand\IH{\operatorname{IH}}

\newcommand{\nc}{\newcommand}
\nc{\renc}{\renewcommand}
\nc{\on}{\operatorname}

\nc\ol{\overline} \nc\wt{\widetilde} \nc\tboxtimes{\wt{\boxtimes}}
\renc{\SS}{{\mathbb S}} \renc{\DD}{{\mathbbD}}
\renewcommand{\AA}{{\mathbb A}}
\nc{\Fq}{{\mathbb F}_q} \nc{\Fqb}{\ol{{\mathbb F}_q}}
\nc{\Ql}{\ol{{\mathbb Q}_\ell}} \renc{\id}{\text{id}}
\nc\X{\mathcal X}

\renc{\Hom}{\on{Hom}} \nc{\Lie}{\on{Lie}} \nc{\Loc}{\on{Loc}}
\nc{\Pic}{\on{Pic}}
\nc{\Sh}{\on{Sh}}
\nc{\pos}{{\on{pos}}} \renc{\Conv}{\on{Conv}} \nc{\Sph}{\on{Sph}}
\renc{\Sym}{\on{Sym}}
\nc{\BunBb}{\overline{\Bun}_B} \nc{\Buno}{\overset{o}{\Bun}}
\nc{\BunPb}{{\overline{\Bun}_P}}
\nc{\BunBM}{\overline{\Bun}_{B(M)}}
\nc{\BunPbw}{{\widetilde{\Bun}_P}}
\nc{\BunBP}{\widetilde{\Bun}_{B,P}} \nc{\GUb}{\overline{G/U}}
\nc{\GUPb}{\overline{G/U(P)}}

\nc{\Hhom}{\underline{\on{Hom}}} \nc\syminfty{\on{Sym}^{\infty}}
\nc\lal{\ol{\lambda}} \nc\xl{\ol{x}} \nc\thl{\ol{\theta}}
\nc\nul{\ol{\nu}} \nc\mul{\ol{\mu}} \nc\Sum\Sigma
\nc{\hl}{\overset{\leftarrow}h}
\nc{\hr}{\overset{\rightarrow}h} \nc{\M}{{\mathcal M}}
\nc{\N}{{\mathcal N}} \nc{\F}{{\mathcal F}} \nc{\D}{{\mathcal D}}
\nc{\Q}{{\mathcal Q}} \nc{\Y}{{\mathcal Y}} \nc{\G}{{\mathcal G}}
\nc{\E}{{\mathcal E}} \nc{\CalC}{{\mathcal C}}
\nc\Dh{\widehat{\D}}
\renewcommand{\O}{{\mathcal O}}
\nc{\C}{{\mathcal C}} \nc{\K}{{\mathcal K}}
\renewcommand{\H}{{\mathcal H}}
\renewcommand{\S}{{\mathcal S}}
\nc{\T}{{\mathcal T}} \nc{\V}{{\mathcal V}} \renc{\P}{{\mathcal
P}} \nc{\A}{{\mathcal A}} \nc{\B}{{\mathcal B}} \nc{\U}{{\mathcal
U}}
\renewcommand{\L}{{\mathcal L}}
\renc{\Gr}{\on{Gr}}

\nc{\frn}{{\check{\mathfrak u}(P)}} \nc{\p}{\mathfrak p}
\nc{\q}{\mathfrak q} \nc\f{{\mathfrak f}}

\nc{\qo}{{\mathfrak q}} \nc{\po}{{\mathfrak p}} \nc{\s}{{\mathfrak
s}} \nc\w{\text{w}}

\nc\mathi\iota \renc\Spec{\on{Spec}} \nc\Mod{\on{Mod}}
\nc{\tw}{\widetilde{\mathfrak t}} \nc{\pw}{\widetilde{\mathfrak
p}} \nc{\qw}{\widetilde{\mathfrak q}} \nc{\jw}{\widetilde j}

\nc{\I}{\mathcal I}

\nc{\lambdach}{{\check\lambda}} \nc{\Lambdach}{{\check\Lambda}{}}
\nc{\much}{{\check\mu}} \nc{\omegach}{{\check\omega}}
\nc{\nuch}{{\check\nu}} \nc{\etach}{{\check\eta}}
\nc{\alphach}{{\check\alpha}} \nc{\betach}{{\check\beta}}
\nc{\rhoch}{{\check\rho}} \nc{\ch}{{\check h}}

\nc{\Hb}{\overline{\H}}

\nc{\BA}{{\mathbb{A}}} \nc{\BC}{{\mathbb{C}}}
\nc{\BG}{{\mathbb{G}}} \nc{\BM}{{\mathbb{M}}}
\nc{\BN}{{\mathbb{N}}} \nc{\BP}{{\mathbb{P}}}
\nc{\BR}{{\mathbb{R}}} \nc{\BZ}{{\mathbb{Z}}}
\nc{\BS}{{\mathbb{S}}}

\nc{\CA}{{\mathcal{A}}} \nc{\CB}{{\mathcal{B}}}

\nc{\CE}{{\mathcal{E}}} \nc{\CF}{{\mathcal{F}}}
\nc{\CG}{{\mathcal{G}}} \nc{\CL}{{\mathcal{L}}}
\nc{\CM}{{\mathcal{M}}} \nc{\CN}{{\mathcal{N}}}
\nc{\CK}{{\mathcal{K}}} \nc{\CO}{{\mathcal{O}}}
\nc{\CP}{{\mathcal{P}}} \nc{\CQ}{{\mathcal{Q}}}
\nc{\CR}{{\mathcal{R}}} \nc{\CS}{{\mathcal{S}}}
\nc{\CT}{{\mathcal{T}}} \nc{\CU}{{\mathcal{U}}}
\nc{\CV}{{\mathcal{V}}} \nc{\CW}{{\mathcal{W}}}
\nc{\CZ}{{\mathcal{Z}}}

\nc{\cM}{{\check{\mathcal M}}{}} \nc{\csM}{{\check{\mathcal A}}{}}
\nc{\obM}{{\overset{\circ}{\mathbf M}}{}}
\nc{\oCA}{{\overset{\circ}{\mathcal A}}{}}
\nc{\obA}{{\overset{\circ}{\mathbf A}}{}}
\nc{\ooM}{{\overset{\circ}{M}}{}}
\nc{\osM}{{\overset{\circ}{\mathsf M}}{}}
\nc{\vM}{{\overset{\bullet}{\mathcal M}}{}}
\nc{\nM}{{\underset{\bullet}{\mathcal M}}{}}
\nc{\obD}{{\overset{\circ}{\mathbf D}}{}}
\nc{\cp}{{\overset{\circ}{\mathbf p}}{}}
\nc{\ofZ}{{\overset{\circ}{\mathfrak Z}}{}}
\nc{\fa}{{\mathfrak{a}}} \nc{\fb}{{\mathfrak{b}}}
\nc{\fg}{{\mathfrak{g}}} \nc{\fgl}{{\mathfrak{gl}}}
\nc{\fh}{{\mathfrak{h}}} \nc{\fj}{{\mathfrak{j}}}
\nc{\fm}{{\mathfrak{m}}} \nc{\fn}{{\mathfrak{n}}}
\nc{\fu}{{\mathfrak{u}}} \nc{\fp}{{\mathfrak{p}}}
\nc{\fr}{{\mathfrak{r}}} \nc{\fs}{{\mathfrak{s}}}
\nc{\fsl}{{\mathfrak{sl}}} \nc{\hsl}{{\widehat{\mathfrak{sl}}}}
\nc{\hgl}{{\widehat{\mathfrak{gl}}}}
\nc{\hg}{{\widehat{\mathfrak{g}}}}
\nc{\chg}{{\widehat{\mathfrak{g}}}{}^\vee}
\nc{\hn}{{\widehat{\mathfrak{n}}}}
\nc{\chn}{{\widehat{\mathfrak{n}}}{}^\vee}

\nc{\fA}{{\mathfrak{A}}} \nc{\fB}{{\mathfrak{B}}}
\nc{\fD}{{\mathfrak{D}}} \nc{\fE}{{\mathfrak{E}}}
\nc{\fF}{{\mathfrak{F}}} \nc{\fG}{{\mathfrak{G}}}
\nc{\fK}{{\mathfrak{K}}} \nc{\fL}{{\mathfrak{L}}}
\nc{\fM}{{\mathfrak{M}}} \nc{\fN}{{\mathfrak{N}}}
\nc{\fP}{{\mathfrak{P}}} \nc{\fU}{{\mathfrak{U}}}
\nc{\fV}{{\mathfrak{V}}} \nc{\fZ}{{\mathfrak{Z}}}

\nc{\bb}{{\mathbf{b}}} \nc{\bc}{{\mathbf{c}}}
\nc{\bd}{{\mathbf{d}}} \nc{\be}{{\mathbf{e}}}
\nc{\bj}{{\mathbf{j}}} \nc{\bn}{{\mathbf{n}}}
\nc{\bp}{{\mathbf{p}}} \nc{\bq}{{\mathbf{q}}}
\nc{\bu}{{\mathbf{u}}} \nc{\bv}{{\mathbf{v}}}
\nc{\bx}{{\mathbf{x}}} \nc{\bs}{{\mathbf{s}}}
\nc{\by}{{\mathbf{y}}} \nc{\bw}{{\mathbf{w}}}
\nc{\bA}{{\mathbf{A}}} \nc{\bK}{{\mathbf{K}}}
\nc{\bB}{{\mathbf{B}}} \nc{\bC}{{\mathbf{C}}}
\nc{\bD}{{\mathbf{D}}} \nc{\bH}{{\mathbf{H}}}
\nc{\bM}{{\mathbf{M}}} \nc{\bN}{{\mathbf{N}}}
\nc{\bV}{{\mathbf{V}}} \nc{\bW}{{\mathbf{W}}}
\nc{\bX}{{\mathbf{X}}} \nc{\bZ}{{\mathbf{Z}}}
\nc{\bS}{{\mathbf{S}}}

\nc{\sA}{{\mathsf{A}}} \nc{\sB}{{\mathsf{B}}}
\nc{\sC}{{\mathsf{C}}} \nc{\sD}{{\mathsf{D}}}
\nc{\sF}{{\mathsf{F}}} \nc{\sK}{{\mathsf{K}}}
\nc{\sM}{{\mathsf{M}}} \nc{\sO}{{\mathsf{O}}}
\nc{\sQ}{{\mathsf{Q}}} \nc{\sP}{{\mathsf{P}}}
\nc{\sZ}{{\mathsf{Z}}} \nc{\sfp}{{\mathsf{p}}}
\nc{\sr}{{\mathsf{r}}} \nc{\sg}{{\mathsf{g}}}
\nc{\sff}{{\mathsf{f}}} \nc{\sfb}{{\mathsf{b}}}
\nc{\sfc}{{\mathsf{c}}} \nc{\sd}{{\mathsf{d}}}

\nc{\BK}{{\bar{K}}}

\nc{\tA}{{\widetilde{\mathbf{A}}}}
\nc{\tB}{{\widetilde{\mathcal{B}}}}
\nc{\tg}{{\widetilde{\mathfrak{g}}}} \nc{\tG}{{\widetilde{G}}}
\nc{\TM}{{\widetilde{\mathbb{M}}}{}}
\nc{\tO}{{\widetilde{\mathsf{O}}}{}}
\nc{\tU}{{\widetilde{\mathfrak{U}}}{}} \nc{\TZ}{{\tilde{Z}}}
\nc{\tx}{{\tilde{x}}} \nc{\tbv}{{\tilde{\bv}}}
\nc{\tfP}{{\widetilde{\mathfrak{P}}}{}} \nc{\tz}{{\tilde{\zeta}}}
\nc{\tmu}{{\tilde{\mu}}}

\nc{\urho}{\underline{\rho}} \nc{\uB}{\underline{B}}
\nc{\uC}{{\underline{\mathbb{C}}}} \nc{\ui}{\underline{i}}
\renc{\uj}{\underline{j}} \nc{\ofP}{{\overline{\mathfrak{P}}}}
\renc{\eps}{\varepsilon} \nc{\hrho}{{\hat{\rho}}}

\nc{\one}{{\mathbf{1}}} \nc{\two}{{\mathbf{t}}}

\nc{\Rep}{{\mathop{\operatorname{\rm Rep}}}}
\nc{\Tot}{{\mathop{\operatorname{\rm Tot}}}}
\renc{\Ker}{{\mathop{\operatorname{\rm Ker}}}}
\nc{\Hilb}{{\mathop{\operatorname{\rm Hilb}}}}
\renc{\End}{{\mathop{\operatorname{\rm End}}}}
\renc{\Ext}{{\mathop{\operatorname{\rm Ext}}}}
\nc{\CHom}{{\mathop{\operatorname{{\mathcal{H}}\it om}}}}
\renc{\GL}{{\mathop{\operatorname{\rm GL}}}}
\nc{\gr}{{\mathop{\operatorname{\rm gr}}}}
\renc{\Id}{{\mathop{\operatorname{\rm Id}}}}
\nc{\de}{{\mathop{\operatorname{\rm def}}}}
\nc{\length}{{\mathop{\operatorname{\rm length}}}}
\nc{\Cliff}{{\mathsf{Cliff}}}
\nc{\Fl}{\on{Fl}} \nc{\Fib}{{\mathsf{Fib}}}
\nc{\Coh}{{\mathsf{Coh}}} \nc{\FCoh}{{\mathsf{FCoh}}}

\nc{\reg}{{\text{\rm reg}}}

\nc{\cplus}{{\mathbf{C}_+}} \nc{\cminus}{{\mathbf{C}_-}}
\nc{\cthree}{{\mathbf{C}_*}} \nc{\Qbar}{{\bar{Q}}}

\nc{\bh}{{\bar{h}}} \nc{\bOmega}{{\overline{\Omega}}}

\nc{\seq}[1]{\stackrel{#1}{\sim}} \nc\QM{{\mathcal {QM}}}

\nc{\chH}{\check H} \nc{\chM}{\check M} \nc{\aff}{{\on{aff}}}

\nc{\chh}{\check \grh}

\renewcommand\chn{\check \grn}

\renewcommand\chg{\check \grg}

\nc\chT{\check T}
\renewcommand\Ext{\operatorname{Ext}}
\newcommand\tBun{\widetilde\Bun}
\newcommand\Eis{\operatorname{Eis}}
\newcommand\Maps{\operatorname{Maps}}
\newcommand\QMaps{\operatorname{QMaps}}

\newcommand{\hinf}{\frac{\infty}{2}}
\newcommand\cht{\check{\grt}}

\title[]{Spaces of quasi-maps into the flag varieties and their applications}
\author[]{Alexander Braverman}

\begin{document}

\begin{abstract}
Given a projective variety $X$ and a smooth projective curve $C$ one may
consider the moduli space of maps $C\to X$. This space admits certain
compactification whose points are called quasi-maps. In the last decade it has been
discovered that in the case when $X$ is a (partial) flag variety of a semi-simple
algebraic group $G$ (or, more generally, of any symmetrizable Kac-Moody Lie algebra)
these compactifications play an important role in such fields as geometric representation theory,
geometric Langlands correspondence, geometry and topology of moduli spaces of $G$-bundles
on algebraic surfaces, 4-dimensional super-symmetric
gauge theory (and probably many others). This paper is a survey of the recent results
about quasi-maps as well as their applications in different branches of representation
theory and algebraic geometry.
\end{abstract}



\section{Introduction}
The spaces of quasi-maps into the flag varieties were introduced by V.~Drinfeld
about 10 years ago and since then proved to play an important role in various
parts of geometric representation theory; more recently it was discovered that
some related constructions are useful also in more classical algebraic geometry
as well as in some questions coming from mathematical physics.

This paper constitutes at attempt to give a more or less self-contained
presentation of the results related to such spaces.
The origin of quasi-maps is as follows: let $C$ be a smooth projective algebraic curve
(over an algebraically closed field $k$) and let $X\subset \PP^N$ be a projective variety
over $k$. One can look at the space $\Maps^d(C,X)$ of maps $C\to X$ such that the composite
map $C\to X\to \PP^N$ has degree $d\in\ZZ_+$. These are quasi-projective schemes of finite type;
in many problems of both representation theory and algebraic geometry it is important to have a
natural compactification of this scheme; one compactification of this sort is provided
by the $\QMaps^d(C,X)$ of {\em quasi-maps} from $C$ to $X$ (cf. Section \ref{def} for the precise
definition). The main property of the scheme $\QMaps^d(C,X)$ is that it possesses a stratification of
the form
\eq{strat}
\QMaps^d(C,X)=\bigcup \limits_{d'=0}^d \Maps^{d'}(C,X)\x \Sym^{d-d'}(C)
\end{equation}
where $\Sym^a(C)$ denote the $a$-th symmetric power of the curve $C$.
In other words, in order to specify a point of $\QMaps^d(C,X)$ one must specify an honest
map $C\to X$ of degree $d'\leq d$ together with $d-d'$ unordered points of $C$.

We must warn the reader from the very beginning that the scheme $\QMaps^d(C,X)$ depends
on the embedding $X\subset \PP^N$. However, in many cases such an embedding is given to us
in the original problem. More generally, when $X$ is a closed subscheme of a product
$\PP^{N_1}\x ...\x \PP^{N_l}$ we may speak about $\QMaps^{d_1,...,d_l}(C,X)$. For example,
if $X$ is the complete flag variety of a semi-simple algebraic group $G$ then $X$ has a canonical
embedding as above (in this case $l$ is the rank of $G$).
We discuss the details in  Section \ref{def}.

We then turn to applications of quasi-maps. In  Section \ref{semi} we explain the relation
between quasi-maps spaces and the so called {\em semi-infinite Schubert varieties}.
In particular, we explain the calculation of the Intersection Cohomology sheaf of
the quasi-maps spaces and relate it to Lusztig's periodic polynomials.
We also mention that quasi-maps could be used to construct some version
of the category of perverse sheaves on the (still not rigorously defined) {\em semi-infinite
flag variety} and relate this category with the category of representations
of the so called {\em small quantum group}.

In Section \ref{eisen} we discuss the results of \cite{BrGa} where the stacks $\oBun_B$
are used in order to construct the so called {\em geometric Eisenstein series} (thus the
contents of \cite{BrGa} have to do with application of the stacks $\oBun_B$
to geometric Langands correspondence; the rest
of the paper is independent of this Section and therefore and can be easily skipped by
a non-interested reader).

In Section \ref{Uhlenbeck} we discuss quasi-maps into affine (partial) flag varieties
and their relation to the Uhlenbeck compactifications of moduli spaces of $G$-bundles on
algebraic surfaces. In Section \ref{nek} we explain how to apply these constructions to certain enumerative questions
related to quantum cohomology of the flag manifolds as well as to N=2 super-symmetric 4-dimensional
gauge theory. Section \ref{open} is devoted to the discussion of some open questions related to the above
subjects.

\smallskip
\noindent
{\bf Acknowledgements.} This paper is mostly based on the author's joint papers  with various people
including
S.~Arkhipov, R.~Bezrukavnikov, P.~Etingof, M.~Finkelberg, D.~Gaitsgory and I.~Mirkovic; I am grateful
to all of them for being very fruitful and patient collaborators. I am also grateful
to J.~Bernstein, V.~Drinfeld and D.~Kazhdan for their constant guidance and to H.~Nakajima, N.~Nekrasov and
A.~Okounkov for interesting and illuminating discussions related to Section \ref{nek} of this paper.

\section{Definition of quasi-maps} \label{def}
In this section we introduce quasi-maps' spaces and some of their relatives.
The reader may skip the details for most applications.

\subsection{Maps and quasi-maps into a projective variety}
Let $X$ be a closed subvariety of the projective space $\PP^N$ and let $C$ be a smooth
projective curve. For any integer $d\geq 0$ we may consider the space $\Maps^d(C,X)$ consisting
of maps $C\to X$ such that the composition $C\to X\to \PP^N$ has degree $d$. This space has
a natural scheme structure and it is in fact quasi-projective. However, it is well-known that
in general it doesn't have to projective (in fact it is almost never projective).

\medskip
\noindent
{\bf Example.} Let $X=\PP^N$. In this case $\Maps^d(C,X)$ classifies the following data:

$\bullet$ A line bundle $\calL$ on $C$ of degree $-d$

$\bullet$ An embedding of vector bundles $\calL\hookrightarrow \calO_C^{N+1}$.

The reason is that every such embedding defines a one-dimensional subspace in $\CC^{N+1}$ for
every point $c\in C$ and thus we get a map $C\to \PP^N$.

Consider, for example, the case when $C=\PP^1$. In that case $\calL$ must be isomorphic
to the line bundle $\calO_{\PP^1}(-d)$ (note that such an isomorphism is defined uniquely
up to a scalar) and thus $\Maps^d(\PP^1,\PP^N)$ becomes an open subset
in the projectivization of the vector space
$\Hom(\calO_{\PP^1}(-d),\calO_{\PP^1}^{N+1})\simeq \CC^{(N+1)(d+1)}$, i.e. $\Maps^d(\PP^1,\PP^N)$ is an open
subset of $\PP^{(N+1)(d+1)-1}$. The reason that it does not coincide with it is that not every non-zero map
$\calO_{\PP^1}(-d)\to \calO_{\PP^1}^{N+1}$ gives rise a map $\PP^1\to \PP^N$ - we need to consider only those
maps which don't vanish in every fiber.

\medskip
\noindent

The above example suggests the following compactification of $\Maps^d(C,X)$. Namely, we define
the space of {\em quasi-maps from $C$ to $X$ of degree $d$} (denoted by $\QMaps^d(C,X)$) to be the
scheme classifying the following data:

1) A line bundle $\calL$ on $C$

2) A non-zero map $\kap: \calL\to \calO_C^{N+1}$

3) Note that $\kap$ defines an honest map $U\to \PP^N$ where $U$ is an open subset of $C$.
We require that the image of this map lies in $X$.

\medskip
\noindent
For example it is easy to see that if $X=\PP^N$ and $C=\PP^1$ then $\QMaps^d(C,X)\simeq \PP^{(N+1)(d+1)-1}$.

In general $\QMaps^d(C,X)$ is projective. Also, set-theoretically it can be explicitly described in the
following way. Assume that we are given a quasi-map $(\calL,\kap)$ as above. Then $\kap$ might
have zeros at points $c_1,...,c_k$ of $C$ of order $a_1,...,a_k$ respectively.
On the other hand, it follows from 3) above that $\kap$ defines an honest map from
the complement to the points $c_1,...,c_k$ to $X$. Since $X$ is projective this map
can be extended to the whole of $C$. Let us call this map $\kap'$. It is easy to see that $\kap'$
has degree $d-\sum a_i$. Also one can recover $\kap$ from $\kap'$ and the collection $(c_1,a_1),...,(c_k,a_k)$.
Thus it follows that $\QMaps^d(C,X)$ is equal to the disjoint union of locally closed subvarieties of the following
form:
\eq{st}
\QMaps^d(C,X)=\bigcup\limits_{0\leq d'\leq d}\Maps^{d'}(C,X)\x \Sym^{d-d'}(C).
\end{equation}
Here $\Sym^{d-d'}(C)$ denotes the corresponding symmetric power of $C$.

Here is a generalization of the above construction. Assume that $X$ is embedded into a product
$\PP^{N_1}\x ...\x\PP^{N_k}$ of projective spaces. Then, in a similar fashion one can talk about
$\Maps^{d_1,...,d_l}(C,X)$ (here all $d_i\geq 0$) and $\QMaps^{d_1,...,d_l}(C,X)$.

\subsection{The case of complete flag varieties}
Let now $G$ be a semi-simple simply connected algebraic group over $k$ and let $\grg$ denote its Lie algebra.
We want to take $X$ to be the complete
flag variety of $G$. If we choose a Borel subgroup $B$ of $G$ then
$X=G/B$. We shall sometimes denote this variety by $X_{G,B}$ (later
we shall also consider the {\em partial flag varieties} $G/P$ associated with
a parabolic subgroup $P\subset G$; this variety will
be denoted by $X_{G,P}$).

Let $V_1,...,V_l$ denote the fundamental representations of $G$. It is well-known that $X_{G,B}$
has a canonical (Pl\"ucker) embedding into $\prod_{i=1}^l \PP(V_i^*)$. This enables
us to talk about quasi-maps into $X$.

We can  describe the set of parameters $(d_1,..,d_l)$ in a little bit
more invariant terms.
First, let us denote by $T$ the Cartan group of $G$ and let $\Lam_G$ denote the coweight
lattice of $G$; by definition $\Lam_G=\Hom(\CC^*,T)$ (in the case $k=\CC$). We have the natural well-known
identification $\Lam_G=H_2(X,\ZZ)$. This allows us to talk about maps $C\to X$ of degree
$\theta\in\Lam_G$. Also if we let $(\ome_1,...,\ome_l)$ denote the
fundamental weights of $G$ then we can also identify
$\Lam_G$ with $\ZZ^l$ by sending a coweight $\theta$ to
$d_1=(\theta,\ome_1),...d_l=(\lam,\theta_l)$. Under these identification $\Maps^{\theta}(C,X)$ is the same
as $\Maps^{d_1,...,d_l}(C,X)$ in the sense of the previous subsection. It is also clear that
this space may be non-empty only if all $d_i\geq 0$. We say that $\theta$ is positive if all $d_i\geq 0$ and denote
the semigroup of all positive $\theta$'s by $\Lam_G^+$.

In the case $C=\PP^1$ we shall denote the space $\Maps^{\theta}(C,X_{G,B})$ by $\calM^{\theta}_{G,B}$ and
the space $\QMaps^{\theta}(C,X_{G,B})$ by $\QM^{\theta}_{G,B}$.

\subsection{Laumon's resolution}
Consider the case $G=SL(n)$ (thus $l=n-1$). In this case $X_{G,B}$ is just the variety of complete flags
$0\subset V_1\subset V_2\subset ...\subset V_{n}=\CC^{n}$, $\dim V_i=i$. Thus, a map $C\to X_{G,B}$ is the same
as a complete flag of subbundles
$$
0\subset\calV_1\subset \calV_2\subset ...\subset \calV_{n}=\calO_C^{n}.
$$
where the rank of $\calV_i$ is equal to $i$. Also we have $d_i=-\deg \calV_i$.

Define now the space $\QMaps^{L,\theta}(C,X_{G,B})$ to consist of all flags as above where $\calV_i$ is an  arbitrary
{\em subsheaf} of $\calO_C^{n}$ of degree $-d_i$. The space was considered by G.~Laumon in \cite{La}.
It is known (cf. \cite{Kuz}) that the natural open embedding of $\Maps^{\theta}(C,X_{G,B})$ into both
$\QMaps^{\theta}(C,X_{G,B})$ and $\QMaps{L,\theta}(C,X_{G,B})$ extends to a projective morphism
$\QMaps^{L,\theta}(C,X_{G,B})\to\QMaps{\theta}(C,X_{G,B})$. In the case $C=\PP^1$ the space
$\QMaps{L,\theta}(\PP^1,X_{G,B})$ is smooth and provides in fact a small resolution of singularities
of $\QMaps^{\theta}(\PP^1,X_{G,B})$.

\subsection{The stacks $\oBun_B$}
Let us fix a curve $C$ as above and let $G$ again be a semi-simple simply connected algebraic
group with a Borel subgroup $B$ (more generally, one can assume that $G$ is any reductive
group whose derived group is simply connected; e.g. one may also consider the
case $G=GL(n)$). We may consider the algebraic stack $\Bun_G=\Bun_G(C)$ classifying
principal algebraic $G$-bundles on $C$. Similarly we may consider the stack
$\Bun_B$ which classifies $B$-bundles. The embedding $B\to G$ gives rise to a the natural
morphism $p:\Bun_B\to \Bun_G$. In the case $G=GL(n)$ the stack $\Bun_G$ classifies vector
bundles of rank $n$ on $C$ and the stack $\Bun_B$ classifies flags of the
forms
$$
0\subset \calV_1\subset\calV_2\subset...\subset\calV_n
$$
where each $\calV_i$ is a vector bundle of rank $i$ on $C$ and the embedding
$\calV_i\to \calV_{i+1}$ are embeddings of vector bundles.

We have the natural projection $B\to T$ (where $T$ as before denotes the Cartan group
of $G$). Hence we also have the natural map $q:\Bun_B\to \Bun_T$. In the case $G=GL(n)$
considered above the group $T$ can be thought of as the group of diagonal matrices;
hence $T$ is naturally isomorphic to $\GG_m^n$.
\footnote{Here $\GG_m$ denotes the multiplicative group}
Thus $\Bun_T$ classifies $n$-tuples $(\calL_1,...,\calL_n)$ of line
bundles on $C$. In terms of the above description of $\Bun_B$ the map
$q$ sends any flag $0\subset \calV_1\subset\calV_2\subset...\subset\calV_n$
to $(\calV_1,\calV_2/\calV_1,...,\calV_n/\calV_{n-1})$.

It is easy to see that in general the connected components of $\Bun_T$ are classified
by elements of the lattice $\Lam_G=\Lam_T$. For each $\theta\in\Lam_G$ we set
$\Bun_B^{\theta}=q^{-1}(\Bun_T^{\theta})$. It is easy to see that the assignment
$\theta\mapsto\Bun_B^{\theta}$ also defines a bijection between
$\Lam_G$ and the set of connected components of $\Bun_B$.

For each $\theta\in\Lam_G$ the map $p:\Bun_B^{\theta}\to\Bun_G$ is representable.
Moreover, it is clear that the fiber of this map over the trivial bundle in $\Bun_G$
is exactly our space $\Maps^{\theta}(C,X_{G,B})$ (note that the stack
$\Bun_B^{\theta}$ exists for any $\theta\in\Lam_G$ but its fiber over the trivial
bundle is non-empty only if $\theta\in\Lam_G^+$). In general, the fibers of $p$
(for fixed $\theta$) are quasi-projective (but not projective) varieties;
for various purposes (discussed, in particular, in other parts of this paper) it is useful
to have a relative compactification $\oBun^{\theta}_B\to\Bun_G$ such that
its fiber over the trivial bundle in $\Bun_G$ will be exactly $\QMaps^{\theta}(C,X_{G,B})$.
Such a compactification indeed can be constructed; let us give its explicit description
(in particular, this will give a slightly different (but equivalent) definition of
$\QMaps^{\theta}(C,X_{G,B})$.

We want to define $\oBun_B$ as a solution to some moduli problem. Since $\oBun_B$ is going to
be an algebraic stack we must define the {\em groupoid} of $S$-points of $\oBun_B$ for any scheme
$S$ over $\CC$.

Let $\Lambdach_G$ be the dual lattice of $\Lam_G$. This is the weight lattice of the group $G$
\footnote{The reader may find this notation a bit bizarre, since usually one uses the
\  $\check{}$-notation for coweights
and not for weights. However, it turns out that  here it is much more convenient
to use our notation; the main reason for this comes from the fact that many results will
be formulated in terms of the Langlands dual group $\check G$ whose {\em weight} lattice is
$\Lam_G$!}
We define
an $S$-point of $\BunBb$ to be a triple $(\F_G,\F_T,\kappa^{\lambdach},\,\,\forall \lambdach\in \Lambdach_G^+)$, where
$\F_G$ and $\F_T$ are as above, and $\kappa^{\lambdach}$ is a map of coherent sheaves
$$
\L_{\F_T}^{\lambdach}\hookrightarrow \V^{\lambdach}_{\F_G},
$$
such that for every geometric point $s\in S$ the restriction $\kappa^{\lambdach}|_{X\times s}$ is an
injection. The last condition is equivalent to saying that $\kappa^{\lambdach}$ is an injection such that the
quotient $\V^{\lambdach}_{\F_G}/\on{Im}(\kappa^{\lambdach})$ is $S$-flat.

The system of embeddings $\kappa^{\lambdach}$ must satisfy  the so-called
{\em Pl\"ucker relations} which can be formulated as follows.

\smallskip

First, for $\lambdach=0$, $\kappa^0$ must be the identity map
$\O\simeq \L^0_{\F_T}\to \V^0_{\F_G}\simeq \O$. Secondly,
for two dominant integral weights $\lambdach$ and $\much$, the map
\begin{equation*}  \label{pluckerone}
\L^{\lambdach}_{\F_T}\otimes\L^{\much}_{\F_T}
\overset{\kappa^{\lambdach}\otimes\kappa^{\much}}\longrightarrow
\V^{\lambdach}_{\F_G}\otimes \V^{\much}_{\F_G}\simeq
(\V^{\lambdach}\otimes \V^{\much})_{\F_G}
\end{equation*}
must coincide with the composition
\begin{equation*} \label{pluckertwo}
\L^{\lambdach}_{\F_T}\otimes\L^{\much}_{\F_T}\simeq \L^{\lambdach+\much}_{\F_T}
\overset{\kappa^{\lambdach+\much}}\longrightarrow
\V^{\lambdach+\much}_{\F_G}\to (\V^{\lambdach}\otimes \V^{\much})_{\F_G}.
\end{equation*}

It is easy to see that if all the maps $\kap^{\lambdach}$ are {\em embeddings of subbundles} (i.e.
$\kap^{\lambdach}$ does not vanish on any fiber over any $c\in C$ then the collection $(\calF_G,\calF_T)$ together
with all
$\kap^{\lambdach}$ defines a point of $\Bun_B$.

Here is another (somewhat more geometric) definition of $\oBun_B$ (note that restricting to the fiber over
the trivial bundle we get yet another definition of $\QMaps^{\theta}(C,X)$.

Let us denote by $U\subset B$ the unipotent radical of $B$. Since we have the natural
isomorphism $G/U=T$ it follows that the variety $G/U$ is endowed with a natural right
action of $T$ (of course, it also has a natural left $G$-action).

It is now easy to  see that the stack $\Bun_B$ classifies the following data:
$$
(\F_G;\F_T;\kappa:\F_G\to G/U{\overset{T}\times}\F_T),
$$
where $\F_G$ is a $G$-bundle, $\F_T$ is a $T$-bundle and $\kappa$ is a $G$-equivariant map.

\smallskip

Recall that $G/U$ is a  quasi-affine variety and let $\GUb$ denote its affine closure.
The groups $G$ and $T$ act on $G/U$ and therefore also on $\GUb$. The basic example of these
varieties that one should keep in mind is the case $G=SL(2)$. In this case
$G/U$ can be naturally identified with $\AA^2\backslash\{ 0\}$ and
$\GUb=\AA^2$ (here $\AA^2$ denotes the affine plane).

\smallskip

We claim now that
an $S$-point of $\BunBb$ is the same as a triple $(\F_G,\F_T,\kappa)$, where
$\F_G$ (resp., $\F_T$) is an $S$-point of $\Bun_G$ (resp., of $\Bun_T$) and
$\kappa$ is a $G$-equivariant map
$$\F_G\to \GUb{\overset{T}\times}\F_T,$$
such that for every geometric point $s\in S$ there is a Zariski-open subset
$C^0\subset C\times s$ such that the map
$$
\kap|_{C^0}:\F_G|_{C^0} \to \GUb{\overset{T}\times}\F_T|_{C^0}
$$
factors through
$G/U{\overset{T}\times}\F_T|_{C^0}\subset \GUb{\overset{T}\times}\F_T|_{C^0}$.

\subsection{Quasi-maps into partial flag varieties}
Let now $P\subset G$ be an arbitrary parabolic subgroup of $G$. Then as before we may consider the
stack $\Bun_P$ of principal $P$-bundles on $C$; this stack is again naturally mapped to $\Bun_G$
and we would like to find some natural relative compactification of it. It turns out that in this case
there exists {\em two different} natural compactifications $\oBun_P$ and $\tBun_P$ such that
the embedding of $\Bun_P$ into both of them extends to a projective morphism
$\tBun_P\to\oBun_P$. We refer the reader to \cite{BrGa} for the corresponding
definitions. Here we shall only explain the geometric source for the existence of
two such compactifications.

As was explained above the stacks $\oBun_B$ are closely related with the varieties
$G/U$ and their affine closures $\GUb$; it is of crucial importance that $G/U$ has
a free $T$-action such that $(G/U)/T=G/B$.

Given a parabolic subgroup $P$ as above one can attach two quasi-affine $G$-varieties
to it: the first one is $G/[P,P]$ and the second one is $G/U_P$ (here $U_P$ denotes the
unipotent radical of $P$; note that if $P$ is a Borel subgroup of $G$
then $[P,P]=U_P$). Let $M$ denote the Levi group of $P$; by definition
$M=P/U_P$. Also, one has the natural isomorphism $P/[P,P]=M/[M,M]$. Thus
the first variety has a natural free action of $M/[M,M]$ and the second has an action
of $M$; moreover, one has $(G/[P,P])/(M/[M,M])=G/P=(G/U_P)/M$. Thus one can use
the quasi-affine closures of $G/[P,P]$ and of $G/U_P$ to construct two
relative compactifications $\oBun_P$ and $\tBun_P$
of the stack $\Bun_P$ in the way similar to what was explained above for $P=B$.

Taking the fibers of the above stacks over the trivial bundle in $\Bun_G$ we get
two different versions of quasi-maps from $C$ to $G/P=X_{G,P}$. In what follows we shall
denote by $\QM^{\theta}_{G,P}$ the space of quasi-maps $\PP^1\to X_{G,P}$ coming from
$\oBun_P$ (the compactification having to do with the variety $G/P$). Here $\theta$ should
be a positive element of the lattice $\Lam_{G,M}$ which is the lattice of
cocharacters of $M/[M,M]$.

\medskip

It turns out that many of the above definitions may be given also when
$\grg$ is replaced by an affine Kac-Moody Lie algebra; the corresponding
spaces of maps and quasi-maps are closely related to  moduli
spaces of $G$-bundles on a rational algebraic surface.
This will be discussed in Section \ref{Uhlenbeck} (for more details the reader should
consult \cite{bfg}).

\section{Quasi-maps into flag varieties and semi-infinite
Schubert varieties}\label{semi}

\subsection{Ordinary Schubert varieties and their singularities}
Let $G$ as a before be a semi-simple simply connected algebraic group and
let $B$ be a Borel subgroup of it.
Recall that we denote $X_{G,B}=G/B$. It is well-known that
the set of $B$-orbits on $X_{G,B}$ is in one-to-one
correspondence with the elements of the Weyl group $W$ of $G$.
For each $w\in W$ we denote the corresponding orbit by $X_{G,B}^w$. It is also known
that each $X_{G,B}^w$ is isomorphic to the affine space $\AA^{\ell(w)}$
where $\ell:W\to\ZZ_+$ is the length function.

The closure $\oX_{G,B}^w\subset X_{G,B}$ of $X_{G,B}^w$ is usually called the Schubert
variety attached to $W$. The singularities of these varieties play a very important
role in various branches of representation theory. It is known (cf. \cite{Fal} and references
therein) that these varieties are normal and have rational singularities.
Let $\IC_{G,B}^w$ denote the intersection cohomology sheaf of $X_{G,B}^w$.
It is also well-known that the stalks of $\IC_{G,B}^w$ can be described in terms of
the {\em Kazhdan-Lusztig polynomials} attached to $W$ (cf. \cite{kl}).

More generally, given two parabolic subgroups $P,Q\subset G$ one may
study the closures of $Q$-orbits on $G/P$; these are the most
general {\em parabolic Schubert varieties}. The stalks of their IC-sheaves
are computed by {\em parabolic Kazhdan-Lusztig polynomials}.

One can generalize the above construction to the loop (or affine) groups
associated with $G$. Namely, given a parabolic subgroup $P$ as above one
may construct two different affine flag varieties $X_{G,P}^{\aff}$ and $\bfX_{G,P}^{\aff}$
associated with the pair $G,P$. We shall refer to the first one as  the
corresponding {\em thin flag affine partial
flag variety} and to the second one
as the {\em thick partial affine flag variety} (in principle, there exist more general partial
affine flag varieties but we shall never consider them in this paper).
Set-theoretically, we can describe $X_{G,P}^{\aff}$ and $\bfX_{G,P}^{\aff}$ as follows.

Let $\calK=\CC((t))$ be the field of formal Laurent power series; let $\calO\subset \calK$
be the ring of Taylor series.
Consider the "loop" group $G(\calK)$ of $\calK$-points of $G$.
Let $I_P\subset G(\calK)$ denote the subgroup of $G(\calO)\subset G(\calK)$
consisting of those Taylor series whose value at
$t=0$ lies in $P\subset G$. When $P=B$ is a Borel subgroup we shall write just $I$ instead of
$I_B$ and call it the Iwahori subgroup of $G((t))$. We shall also denote by
$I^0\subset I$ its pro-unipotent radical (it consists of those Taylor series as above whose value
at $t=0$ lies in the unipotent radical $U$ of $B$). Note also that when
$P=G$ we have $I_G=G(\calO)$.

Similarly, we can define the group $\bfI_P\subset G[t^{-1}]$ consisting of those polynomials
in $t^{-1}$ whose value at $t=\infty$ lies in $P$. Thus on level of $\CC$-points
we have
$$
X_{G,P}^{\aff}=G(\calK)/I_P;\quad \bfX_{G,P}^{\aff}=G(\calK)/\bfI_P.
$$

We shall be mostly talking talking about Schubert varieties in $X_{G,P}^{\aff}$
(the flag varieties $\bfX_{G,P}^{\aff}$ will also appear in Section \ref{Uhlenbeck} of this
paper). By definition, these are closures of $I_Q$-orbits in some $X_{G,P}^{\aff}$.
These are known to be finite-dimensional normal projective varieties having
rational singularities (cf. \cite{Fal}). In the case when $P=Q=B$ (i.e. when we are dealing
with $I$-orbits on $X_{G,B}^{\aff}$) the orbits are classified by
elements of the affine Weyl group $W_{\aff}$ (by definition, this is the semi-direct product
of the Weyl group $W$ of $G$ and the lattice $\Lam_G$). At the other extreme, when
$P=Q=G$ we are dealing with $G(\calO)$-orbits on $G(\calK)/G(\calO)$; these orbits are in
one-to-one correspondence with $\Lam_G/W$. The latter set can be identified with the set
of dominant weights of the Langlands dual group ${\check G}$.

One of the reasons that the complete flag variety $X_{G,B}$ plays a distinguished role in
representation theory is the Beilinson-Bernstein localization theorem (cf. \cite{bb}) which allows one to realize
representations of the Lie algebra $\grg$ (with fixed central character of the universal
enveloping algebra $U(\grg)$) in terms of algebraic $D$-modules on $X_{G,B}$.
In this way the category of $B$-equivariant (or, more generally, $U$-equivariant) modules
corresponds to the regular block of the so called category $\calO$. Similar (but much less
understood) statements hold in the affine case too. Namely let $\grg_{\aff}$ denote the
affine Lie algebra corresponding to $\grg$. By definition this algebra is a central extension
of the loop  algebra $\grg((t))$:
$$
0\to \CC\to \grg_{\aff}\to \grg((t))\to 0.
$$
Then one gets a geometric realization
(of some part of) the category for the corresponding affine Lie algebra $\grg_{\aff}$;
the variety $X_{G,B}^{\aff}$ allows to realize $\grg_{\aff}$-modules at the negative level
and the variety $\bfX_{G,B}^{\aff}$ has to do with $\grg_{\aff}$-modules on the positive level.
We refer the reader to \cite{Kash} and references therein for more details.

\subsection{The semi-infinite flag manifold}

The semi-infinite flag manifold $X^{\hinf}_{G,B}$ is usually defined as the quotient
$G(\calK)/T(\calO)\cdot U(\calK)$. In terms of algebraic geometry this "space" seems to
be widely infinite-dimensional. However, one still would like to think of it
as some kind of geometric object; this should have many applications to representation theory.

More specifically, we would like to mention the following two problems:

1) Construct the category of $D$-modules (or perverse sheaves) on $X_{G,B}^{\hinf}$
and relate it to some other abelian categories coming from representation theory
of affine Lie algebras and quantum groups.

\smallskip

2) It is easy to see that the orbits of the Iwahori group $I$ on $X_{G,B}^{\hinf}$ are classified
by elements of the affine Weyl group $W_{\aff}$ attached to $G$. We shall denote by
$X_{G,B}^{w,\hinf}$ the orbit corresponding to $w\in W_\aff$. Then the problem reads as follows:
explain in what sense the singularities of the closures of $\oX_{G,B}^{w,\hinf}$ are finite-dimensional
and "understand" those singularities. In particular, one should be able to compute the stalks
of the IC-sheaves associated to those singularities and relate them to the periodic
polynomials defined in \cite{Lper}. More generally, one can study $I_P$-orbits on $X_{G,B}^{\hinf}$
together with their closures. We shall refer to them as the (not yet constructed) {\em semi-infinite
Schubert varieties}.

In particular, the $G(\calO)$-orbits on $X_{G,B}^{\hinf}$ are
classified by elements of $\Lam_G$; for each $\mu\in\Lam_G$ we shall denote the corresponding
orbit simply by $S^{\mu}$. It is easy to see that if the closure $\oS^{\mu}$ of $S^{\mu}$ makes any
reasonable sense then it must be equal to the union of all $S^{\nu}$'s with
$\nu-\mu\in\Lam_G^+$. Also, the lattice $\Lam_G=T(\calK)/T(\calO)$ acts on $X_{G,B}^{\hinf}$
on the right and every $\gam\in \Lam_G$ maps the orbit $S^{\mu}$ to $S^{\mu+\gam}$.
Hence the singularity of $\oS^{\mu}$ in the neighborhood
of a point of $S^{\nu}$ depends only on $\nu-\mu\in\Lam_G^+$. In particular, if the intersection cohomology
sheaf $\IC(\oS^{\mu})$ makes sense, then its stalk at a point of $S^{\nu}$ should only depend
on $\nu-\mu=\theta$. In fact, from the results of \cite{Lqan}, \cite{Lper} and references therein it is natural
to expect that this stalk comes from the graded vector space $U_{\theta}$ computed as follows:
let $\chg$ denote the Langlands dual Lie algebra whose root system is dual to that of $\grg$.
We have its triangular decomposition $\chg=\chn_-\oplus\cht\oplus \chn_+$.
Also $\cht=\grt^*$ and  we may identify $\Lam_G$ with the root lattice $\chg$.
Consider the symmetric algebra $\Sym(\grn_+)$ with the natural even grading on it (defined
by the requiring that the subspace $\chn_+\subset \Sym(\grn_+)$ has degree 2).
The dual Cartan torus ${\check T}$ acts on this algebra (since it acts on $\grn_+$);
for each $\theta\in\Lam_G^+$ we may consider the subspace $\Sym(\grn_+)_\theta\subset \Sym(\grn_+)$
on which ${\check T}$ acts by the character $\theta$. This space inherits the grading from
$\Sym(\grn_+)$. Let also $\check{\rho}\in\check{\Lam}_G$ denote the half-sum of the positive roots of
$G$. Then, guided by the results
of {\em loc. cit.} one expects to have
\eq{utheta}
U_{\theta}\simeq \Sym(\grn_+)_{\theta}[2\la\theta,\check{\rho}\ra].
\end{equation}

\medskip
The general principle (due to Drinfeld) says that one should be able to use the quasi-maps spaces $\QM_{G,B}^{\theta}$
(or, the stacks $\oBun_B$ for any smooth projective curve $C$) as "finite-dimensional models" for the semi-infinite
flag manifolds and the semi-infinite Schubert varieties. In particular, one expects to be able construct  the
correct category of $D$-modules using quasi-maps as well as to turn \refe{utheta} into a mathematical theorem.
This has indeed been performed in the works \cite{ABBGM}, \cite{FFKM} and \cite{BFGM}.
Let us give a brief sketch of the results of {\em loc. cit.}

\subsection{Localization theorem for the small quantum group}
Let us turn to some illustrations of the above principle. First of all,
in \cite{ABBGM} we propose a definition of the category
$\Perv(X_{G,B}^{\hinf})$ in terms of the stacks $\oBun_B$. We give a representation-theoretic
interpretation of the corresponding subcategory $\Perv_{I^0}(X_{G,B}^{\hinf})$ consisting
of $I^0$-equivariant perverse sheaves; it turns out to be equivalent to the {\em regular block} of
category of graded representations of the so called {\em small quantum group} $\gru_\ell$ attached
to the Lie algebra $\grg$; here $\ell$ denotes a root of unity satisfying some mild assumptions
(cf. \cite{ABBGM} for more details). This result was conjectured by B.~Feigin and E.~Frenkel in the early 90s.
Another representation-theoretic interpretation of the same category (in terms of representations
of the affine Lie algebra $\grg_{\aff}$) should appear soon in the works of Frenkel and Gaitsgory.

\subsection{Computation of the IC-sheaf}
Another check of the above principle will be to compute the stalks of the IC-sheaves of the spaces
$\QM^{\theta}_{G,B}$ (or the stacks $\oBun_B$) and compare it with \refe{utheta} .
This was done in \cite{FFKM}; also in \cite{BFGM} this was generalized to arbitrary
parabolic $P\subset G$. More specifically, the space $\QM^{\theta}_{G,B}$ possesses the following
stratification (similar to \refe{st}):
\eq{st-fl}
\QM^{\theta}_{G,B}=\bigcup\limits_{\mu\in\Lam_G^+}\calM^{\theta-\mu}_{G,B}\x \Sym^{\mu}(\PP^1).
\end{equation}
Here by $\Sym^{\mu}(\PP^1)$ we mean the space of all colored divisors $\sum \mu_i x_i$ where
$\mu_i\in\Lam_G^+, x_i\in\PP^1$ and $\sum \mu_i=\mu$. Then we have
\begin{theorem}\label{IC}
The stalk of $\IC_{\QM^{\theta}_{G,B}}$ at a point of $\QM^{\theta}_{G,B}$ corresponding to
a colored divisor $\sum \mu_i x_i$ as above is equal to
$\otimes_i \Sym(\chn_+)_{\mu_i}[2\la\mu_i,\check{\rho}\ra]$. In particular, the
stalk of $\IC_{\QM^{\theta}_{G,B}}$ at the "most singular" points of $\QM^{\theta}_{G,B}$
corresponding to the divisor of the form $\theta\cdot x$ (for some $x\in \PP^1$) is equal to
$U_{\theta}=\Sym(\grn_+)_{\theta}[2\la\theta,\check{\rho}\ra]$.
\end{theorem}
The proof of Theorem \ref{IC} relies on many things, in particular the results of \cite{MV}
about semi-infinite orbits in the affine Grassmannian of $G$.

\subsection{Geometric construction of the universal Verma module}\label{g-verma}
We have seen that one can read off some information
related to the Langlands dual Lie algebra from the singularities of
the quasi-maps' spaces. It is natural to ask
if one could push this a little further and get a geometric
construction of $\chg$-modules (in Section \ref{Uhlenbeck} we are going to generalize
it to affine Lie algebras).

Of course, the most interesting modules that one would like to get
in this way are the finite-dimensional modules.
This, however, has not been done yet.
In this section we explain how to use the spaces of quasi-maps
in order to construct the "universal
Verma module" for the Lie algebra $\chg$. We also
give geometric interpretation of the Shapovalov form and the Whittaker
vectors (cf. the definitions below). We shall generalize this
in Subsection \ref{uhl-verma} to the case of affine Lie algebras. These constructions
will play the crucial role in Section \ref{nek} where we discuss
applications of our techniques to some questions of enumerative algebraic
geometry.

First, let $Y$ be a scheme endowed with an action of a reductive algebraic group $L$ (in most applications
$L$ will actually be a torus). We denote by $\IH_L(Y)$ the intersection cohomology of $Y$ with complex coefficients.
This is a module over the algebra $\calA_{L}=H^*_L(pt)$ which is known to be isomorphic to the algebra of polynomial
functions on the Lie algebra $\grl$ of $L$ which are invariant under the adjoint action of $L$.
We let $\calK_L$ denote the field of fractions of $\calA_L$.

We now take $Y$ to be the space $^{\bt}\QM_{G,B}^{\theta}$ of {\em based} quasi-maps $\PP^1\to X_{G,B}$. By definition,
this is the locally closed subscheme of $\QM^{\theta}_{G,B}$ corresponding to those quasi-maps which
are first of all well-defined as {\em maps} around $\infty\in\PP^1$ and such that their value at $\infty$ is
equal to the point $e_{G,B}\in X_{G,B}$ corresponding to the unit element of $G$ under the identification
$X_{G,B}=G/B$. This scheme is endowed with a natural action of the torus $T\x \CC^*$ (here
$T$ acts on $X_{G,B}$ preserving $e_{G,B}$ and $\CC^*$ acts on $\PP^1$ preserving $\infty$). Define
$$
\IH_{G,B}^{\theta}=\IH^*_{T\x\CC^*}
(^{\bt}\QM_{G,B}^{\theta})\underset{\calA_{T\x\CC^*}}
\ten \calK_{T\x\CC^*}, \quad
\IH_{G,B}
=\bigoplus\limits_{\theta\in \Lam_{G,B}^+}\IH_{G,B}^{\theta}
$$
Each $\IH_{G,B}^{\theta}$ is a finite-dimensional vector space over the field
$\calK_{T\x\CC^*}$ which can be thought of as the field of rational functions
of the variables $a\in\grt$ and $\hbar\in \CC$. Moreover, $\IH_{G,B}^{\theta}$ is
endowed with a (non-degenerate) Poincar\'e pairing $\la\cdot,\cdot\ra_{G,P}^{\theta}$
taking values in
$\calK_{T\x\CC^*}$ (one has to explain why the Poincar\'e pairing is well defined since $^{\bt}\QM_{G,B}^{\theta}$
is not projective; this is a corollary of (some version of) the localization theorem in equivariant
cohomology - cf. \cite{Br} for more details).

In \cite{Br} we construct a natural action of the Lie algebra $\chg$
on the space $\IH_{G,B}$. Moreover, this action has the following
properties. First of all, let us denote by
$\la\cdot,\cdot\ra_{G,B}$ the direct sum of the pairings
$(-1)^{\la \theta,\check{\rho}\ra}\la\cdot,\cdot\ra_{G,B}^{\theta}$.

Recall that the Lie algebra $\chg$ has its triangular
decomposition
$\chg=\chn_+\oplus\cht\oplus\chn_-$. Let $\kap:\chg\to\chg$ denote the
Cartan anti-involution which interchanges $\chn_+$ and $\chn_-$ and
acts as identity on $\cht$.
For each $\lam\in\grt=(\cht)^*$ we denote by $M(\lam)$ the
corresponding
Verma module with lowest weight $\lam$; this is a module generated by
a vector $v_\lam$ with (the only) relations
$$
t(v_\lam)=\lam(t)v_\lam\quad\text{for $t\in\cht$ and}\quad
n(v_\lam)=0\quad\text{for $n\in\chn_-$}.
$$

\noindent
Then:

1) $\IH_{G,B}$ (with the
above action) becomes isomorphic to $M(\lam)$ where
$\lam=\frac{a}{\hbar}+\rho$.

2) $\IH_{G,B}^{\theta}\subset \IH_{G,B}$ is the
$\frac{a}{\hbar}+\rho+\theta$-weight space of $\IH_{G,B}$.

3) For each $g\in\chg$ and $v,w\in \IH_{G,B}$ we have
$$
\la g(v),w\ra_{G,B}=\la v,\kap(g)w\ra_{G,B}.
$$

4) The  vector
$\sum_{\theta}1_{G,B}^{\theta}$ (lying in  the corresponding completion of
$\IH_{G,B}$) is a Whittaker vector (i.e. an $\grn_-$-eigen-vector) for the above action.

We are not going to explain the construction of the action in this survey
paper. Let us only make a few remarks about it. In the case $G=SL(n)$ the smallness
result of \cite{Kuz} allows to
replace the intersection cohomology of $^{\bt}\QM_{G,B}^{\theta}$ by the ordinary
cohomology of the corresponding based version of the Laumon resolution $^{\bt}\QM_{G,B}^{L,\theta}$;
on the latter (equivariant localized) cohomology the action of the Chevalley generators of $\chg=sl(n)$ can be defined
by means of some explicit correspondences (this is similar to the main construction of
\cite{fk}; also in \cite{bf} we generalize this to the case when equivariant cohomology is replaced by
equivariant $K$-theory. In this case the action of the Lie algebra $sl(n)$ is replaced by the action
of the corresponding quantum group $U_q(sl(n))$).
Also for any $G$ the fact, that the dimension of $\IH_{G,B}^{\theta}$ can be easily
deduced from Theorem \ref{IC}.
Our  construction of the $\chg$-action on $\IH_{G,B}$ is very close
to the construction in Section 4 of \cite{FFKM}.

\section{The stack $\oBun_B$ and geometric Eisenstein series}\label{eisen}
\label{eisen}
This section is devoted to an application of the stacks $\oBun_B$  to some questions of {\em geometric
Langlands correspondence}. A reader who is not interested in the
subject may skip this section since it will never be used in the
future. In fact we are going to discuss only one such application (which was the first one
historically) - the construction of {\em geometric Eisenstein series}.
Let us note, though, that
the stacks $\oBun_B$ have appeared in many other works on the subject.
For example they play the  crucial role in the geometric proof of Casselman-Shalika
formula by E.~Frenkel, D.~Gaitsgory and K.~Vilonen (cf. \cite{FGV}) as well as in the proof of
 the so called "vanishing conjecture" which implies (the main portion of) the geometric
Langlands conjecture for $GL(n)$ (cf. \cite{FGV1}, \cite{Ga}) and the de Jong conjecture about representations
of Galois groups of functional fields (cf. \cite{GadJ}). A good review of these results may
be found in \cite{GaICM}.

All the results discussed below are taken from
\cite{BrGa}.

\subsection{The usual Eisenstein series}
Let $X$ be a curve over $\Fq$ and let $G$ be a reductive group. The classical
theory of automorphic forms is concerned with the space of functions on
the quotient $G_{\AA}/G_\calK$, where
$\calK$ (resp., $\AA$) is the field or rational functions on $X$
(resp., the ring of ad\`eles of $\calK$).
In this paper, we will consider only the unramified situation, i.e.
we will study functions
(and afterwards perverse sheaves) on the double quotient
$G_{\OO}\backslash G_{\AA}/G_\calK$.

Let $T$ be a Cartan subgroup of $G$. There is a well-known construction,
called the {\it Eisenstein series operator}
that attaches to a compactly supported
function on $T_{\OO}\backslash T_{\AA}/T_\calK$ a function on
$G_{\OO}\backslash G_{\AA}/G_\calK$:

Consider the diagram
$$
\begin{CD}
B_{\OO}\backslash B_{\AA}/B_\calK   @>{q}>> T_{\OO}\backslash T_{\AA}/T_\calK \\
@V{p}VV   \\
G_{\OO}\backslash G_{\AA}/G_\calK,
\end{CD}
$$
where $B$ is a Borel subgroup of $G$. Up to a normalization factor,
the Eisenstein series of a function $S$ on
$T_{\OO}\backslash T_{\AA}/T_\calK$ is ${p}_{!}({q}^*(S))$, where
${q}^*$ denotes pull-back and ${p}_{!}$ is integration along the fiber.

Our goal is to study a geometric analog of this construction.

Let $\Bun_G$ denote the stack of $G$--bundles on $X$. One may regard the derived category
of constructible sheaves on $\Bun_G$ (denoted $\Sh(\Bun_G)$) as a
geometric analog of the space of functions on $G_{\OO}\backslash G_{\AA}/G_\calK$. Then, by geometrizing
the Eisenstein series operator, we obtain an Eisenstein series functor
$\Eis'$%:\Sh(\Bun_T)\to\Sh(\Bun_G)$ defined by a diagram
similar to the above one, where the intermediate stack is $\Bun_B$--the stack of $B$--bundles on $X$.

However, this construction has an immediate drawback--it is not
sufficienly functorial (for example
it does not commute with Verdier duality),
the reason being that the projection
${p}:\Bun_B\to\Bun_G$ has non-compact fibers.
Therefore, it is natural to look for a relative
compactification of $\Bun_B$ along the fibers of the
projection ${p}$.

It turns out that the compactification  $\BunBb$ discussed indeed does the job, i.e. we can use it to define
the corrected functor $\Eis:\Sh(\Bun_T)\to\Sh(\Bun_G)$. The paper \cite{BrGa} is devoted to the
investigation of various properties of this functor.

In fact, all the technical results about the functor $\Eis$ essentially reduce to questions about the geometry
of $\BunBb$ and the behaviour of the intersection cohomology sheaf on it.

We should say right away that the pioneering work in this direction was done by G.~Laumon in \cite{La},
who considered the case of $G=GL(n)$ using his own compactification  $\oBun_B^L$ of the stack $\Bun_B$.
In the sequel we will explain how the two approaches are related.

\subsection{Survey of the main results of \cite{BrGa}}

Once the stack $\BunBb$ is constructed, one can try to use it to define the "compactified" Eisenstein series
functor $\Eis:\Sh(\Bun_T)\to\Sh(\Bun_G)$. Let $p$ and $q$ denote the natural projections from $\BunBb$
to $\Bun_G$ and $\Bun_T$, respectively. The first idea would be to consider the functor
$\S\in\Sh(\Bun_T)\mapsto p_{!}(q^*(\S))\in\Bun_G$. However, this is too naive, since if we want
our functor to commute with Verdier duality, we need to take into account the singularities of $\BunBb$.
Therefore, one introduces a {\it kernel} on $\BunBb$ given by its intersection cohomology sheaf. I.e., we define
the functor $\Eis$ by
$$\S\mapsto p_{!}(q^*(\S)\otimes\IC_{\BunBb}),$$
up to a cohomological shift and Tate's twist.
Similarly, one defines the functor $\Eis^G_M:\Sh(\Bun_M)\to\Sh(\Bun_G)$, where $M$ is the Levi quotient
of a parabolic $P$.

The first test whether our definition of the functor $\Eis$ is "the right one" would be the assertion
that $\Eis$ (or more generally $\Eis^G_M$) indeed commutes with Verdier duality.
It can be shown that our $\Eis$ indeed passes this test.

Let us again add a comment of how the functor $\Eis$ is conneceted to Laumon's work. One can define
functors $\Eis^L:\Sh(\Bun_T)\to\Sh(\Bun_G)$ using Laumon's compactification. (In the original work \cite{La},
Laumon did not consider $\Eis^L$ as a functor, but rather applied it to specific sheaves on $\Bun_T$.)
However, from the smallness result of \cite{Kuz} it follows that the functors $\Eis^L$ and $\Eis$ are canonically
isomorphic.

Once we defined the functors $\Eis=\Eis^G_T:\Sh(\Bun_T)\to\Sh(\Bun_G)$, $\Eis^G_M:\Sh(\Bun_M)\to\Sh(\Bun_G)$
and a similar functor for $M$, $\Eis^M_T:\Sh(\Bun_T)\to\Sh(\Bun_M)$, it is by all means natural to expect
that these functors compose nicely, i.e. that $\Eis^G_T\simeq \Eis^G_M\circ \Eis^M_T$.

For example, if instead of $\Eis$ we used the naive (uncompactified)
functor $\Eis'$, the analogous assertion would be a triviality, since
$\Bun_B\simeq \Bun_P\underset{\Bun_M}\times\Bun_{B(M)}$, where $B(M)$ is the Borel subgroup of $M$.

The problem with our definition of $\Eis^G_M$ is that there is no map between the relevant
compactifications, i.e. from $\BunBb$ to $\BunPbw$. Neverthess, the assertion that
$\Eis^G_T\simeq \Eis^G_M\circ\Eis^M_T$ does hold. This in fact is a non-trivial theorem proved in \cite{BrGa}.

Here are the main properties of the Eisenstein series functor.

\subsection{Behaviour with respect to the Hecke functors}
Classically, on the space of functions on the double quotient $G_{\OO}\backslash G_{\AA}/G_\calK$
we have the action of $\underset{x\in X}\otimes \H_x(G)$, where $x$ runs over the set of places of $\calK$,
and for each $x\in X$, $\H_x(G)$ denotes the corresponding spherical Hecke algebra of the group $G$.

Similarly, $\underset{x\in X}\otimes \H_x(T)$ acts on the space of functions on
$T_{\OO}\backslash T_{\AA}/T_\calK$. In addition, for every $x$ as above, there is a canonical homomorphism
$\H_x(G)\to \H_x(T)$ described as follows:

Recall that there is a canonical isomorphism (due to Satake) between $\H_x(G)$ and the Grothendieck
ring of the category of finite-dimensional representations of the Langlands dual group $\check G$.
We have the natural restriction functor $\on{Rep}(\check G)\to \on{Rep}(\check T)$, and our homomorphism
$\H_x(G)\to \H_x(T)$ corresponds to the induced homomorphism $K(\on{Rep}(\check G))\to K(\on{Rep}(\check T))$
between Grothendieck rings.

The basic property of the Eisenstein series operators is that it intertwines the $\H_x(G)$-action on
$G_{\OO}\backslash G_{\AA}/G_\calK$ and the $\H_x(T)$-action on $T_{\OO}\backslash T_{\AA}/T_\calK$
via the above homomorphism.

\smallskip

Our result below  is a reflection of this phenomenon in the geometric setting. Now, instead
of the Hecke algebras, we have the action of the Hecke functors on $\Sh(\Bun_G)$. Namely,
for $x\in X$ and an object $V\in \on{Rep}(\check G)$, one defines the {\it Hecke functor}
$$
\S\mapsto {}_xH_G(V,\S)
$$
from $\Sh(\Bun_G)$ to
itself. The existence of such functors comes from the so called {\em geometric Satake isomorphism} -
cf. \cite{MV} and references therein.

We claim that for any $\S\in\Sh(\Bun_T)$ we have:
$$
_xH_G(V,\Eis(\S))\simeq \Eis({}_xH_T(\on{Res}^G_T(V),\S)).
$$
This theorem is more or less equivalent to one of the main results of \cite{MV}.
A similar statement holds for the non-principal Eisenstein series functor $\Eis^G_M$.

As a corollary, we obtain that if $E_{\check M}$ is an $\check M$-local system on $X$ and
$\Aut_{E_{\check M}}$ is a perverse sheaf (or a complex of sheaves) on $\Bun_M$,
corresponding to it in the sense of the geometric Langlands correspondence, then the complex
$\Eis^G_M(\Aut_{E_{\check M}})$ on $\Bun_G$ is a Hecke eigen-sheaf with respect to the induced
$\check G$-local system.

In particular, we construct Hecke eigen-sheaves
for those homomorphisms $\pi_1(X)\to \check G$, whose image is contained in a maximal torus of $\check G$.

\subsection{The functional equation}\label{func}

It is well-known that the classical Eisenstein series satisfy the functional equation. Namely, let $\chi$
be a character of the group $T_{\OO}\backslash T_{\AA}/T_\calK$ and let $\w\in W$ be an element of the Weyl
group. We can translate $\chi$ by menas of $\w$ and obtain a new (grossen)-character $\chi^{\w}$.

The functional equation is the assertion that the Eisenstein series corresponding to $\chi$ and
$\chi^{\w}$ are equal, up to a ratio of the corresponding L-functions.

Now let $\S$ be an arbitrary complex of sheaves on $\Bun_T$ and let $\w\cdot\S$ be its $\w$-translate.
One may wonder whether there is any relation between $\Eis(\S)$ and $\Eis(\w\cdot\S)$.

We single out a subcategory in $\Sh(\Bun_T)$, corresponding to sheaves which we call "regular", for
which we answer the above question. We show that for a regular sheaf $\S$ we have
$$
\Eis(\S)\simeq \Eis(\w\cdot\S).
$$

(It is easy to see that one should not expect the functional equation to hold for
non-regular sheaves.)

A remarkable feature of this assertion is that the L-factors that enter the classical functional
equation have disappeared. An explanation of this fact is provided by the corresponding result
from \cite{BrGa}
which says that the definition of $\Eis$ via the intersection cohomology sheaf on $\BunBb$ already
incorporates the $L$--function.

We remark that an assertion similar to the above functional equation should hold also for
non-principal Eisenstein series. Unfortunately, this seems to be beyond the access of our methods.

Using the above results we obtain a proof
of the following very special case of the Langlands conjecture.  Namely, we prove that if
we start with an unramified irreducible representation of $\pi_1(X)$ into $\check G$,
such that $\pi_1(X)^{geom}$
\footnote{Here $\pi_1(X)^{geom}$ denotes the geometric fundamental group of $X$ - i.e. the fundamental group
of $X$ over the algebraic closure of $\FF_q$}
maps to $\check T\subset \check G$, then there
exists an unramified automorphic form on $G_{\AA}$ which corresponds to this representation in the
sense of Langlands.

This may be considered as an application of the machinery
developed in \cite{BrGa} to the classical theory of automorphic forms.

\section{Quasi-maps into affine flag varieties and  Uhlenbeck
compactifications}\label{Uhlenbeck}

In this Section we take the base field to be $\CC$.

\subsection{The problem}
Let $G$ be an almost simple simply connected group over $\BC$, with
Lie algebra $\fg$, and let $\bS$ be a smooth projective surface.

Let us denote by $\Bun^d_G(\bS)$ the moduli space (stack) of principal
$G$-bundles on $\bS$ of second Chern class $d\in\ZZ$. It is easy to see
that $\Bun^d_G(\bS)$ cannot be compact and for many reasons it is natural to expect
that there exists a  compactification of
$\Bun^d_G(\bS)$ which looks like  a union
\eq{Uhlenbeck?}
\underset{b\in \BN}\bigcup\, \Bun^{d-b}_G(\bS)\times \on{Sym}^b(\bS).
\end{equation}
Note the striking similarity between \refe{Uhlenbeck?} and \refe{st}.

In the differential-geometric framework of moduli spaces of $K$-instantons
on Riemannian 4-manifolds (where $K$ is the maximal compact subgroup of $G$)
such a compactification was introduced in the pioneering work ~\cite{u}.
Therefore, we shall call its algebro-geometric version the
Uhlenbeck space, and denote it by $\calU^d_G(\bS)$.
%S.~Donaldson ~\cite{d},
%following the pioneering work ~\cite{u}.
%In the framework of algebraic geometry,
%Such compactifications have been constructed in differential geometry
%by K.~Uhlenbeck for $K$-instantons (where $K$ is the maximal compact
%subgroup of $G$), cf. \cite{u}, \cite{d}.

\medskip

Unfortunately, one still does not know how to construct the spaces
$\calU^d_G(\bS)$ for a general group $G$ and an arbitrary surface $\bS$.
More precisely, one would like to formulate a moduli problem, to
which $\calU^d_G(\bS)$ would be the answer, and so far this is not known.
In this formulation the question of constructing the Uhlenbeck spaces
has been posed (to the best of our knowledge) by V.~Ginzburg. He and
V.~Baranovsky (cf. \cite{BaGi}) have made the first attempts to
solve it, as well as indicated the approach adopted in this paper.

\medskip

A significant simplification occurs for $G=SL_n$.
Let us note that when $G=SL_n$, there exists another natural
compactification of the stack $\Bun^d_n(\bS):=\Bun^d_{SL_n}(\bS)$
(called the Gieseker compactification),
by torsion-free sheaves of generic rank $n$ and of second Chern
class $a$, called the Gieseker compactification, which in this paper
we will denote by $\wt{\fN}^d_n(\bS)$. One expects that there exists
a proper map $\sff:\wt{\fN}^{d}_n(\bS)\to \calU^{d}_{SL_n}(\bS)$,
described as follows:

A torsion-free sheaf $\M$ embeds into a short exact sequence
$$0\to \M\to \M'\to \M_0\to 0,$$
where $\M'$ is a vector bundle (called the saturation of $\M$), and
$\M_0$ is a finite-length sheaf. The map $\sf$ should send a point
of $\wt{\fN}^d_n(\bS)$ corresponding to $\M$ to the pair
$(\M',\on{cycle}(\M_0))\in \Bun^{d-b}_n(\bS)\times \on{Sym}^b(\bS)$,
where $b$ is the length of $\M_0$, and $\on{cycle}(\M_0)$ is the cycle
of $\M_0$. In other words, the map $\sf$ must ``collapse" the information
of the quotient $\M'\to \M_0$ to just the information of the length of
$\M_0$ at various points of $\bS$.

Since the spaces $\wt{\fN}^d_n(\bS)$, being a solution of a moduli problem,
are easy to construct, one may attempt to construct the Uhlenbeck spaces
$\calU^d_{SL_n}(\bS)$ by constructing an explicit blow down of the
Gieseker spaces $\wt{\fN}^d_n(\bS)$. This has indeed been performed in
the works of J.~Li (cf. \cite{Li}) and J.~W.~ Morgan (cf. \cite{Mo}).

The problem simplifies even further, when we put $\bS=\BP^2$, the
projective plane, and consider bundles trivialized along a fixed
line $\BP^1\subset \BP^2$.
In this case, the sought-for space $\calU^d_n(\bS)$
has been constructed by S.~Donaldson  and thoroughly
studied by H.~Nakajima (cf. e.g. \cite{n1}) in his works on quiver
varieties.

\medskip

In \cite{bfg} we  consider the case of an arbitrary group
$G$, but the surface equal to $\BP^2$ (and we will be interested in
bundles trivialized along $\BP^1\subset \BP^2$, i.e., we will work
in the Donaldson-Nakajima set-up.)

In fact we are able to construct the Uhlenbeck spaces $\calU^d_G$, but only
up to nilpotents. I.e., we will have several definitions, two of
which admit modular descriptions, and which produce the same answer
on the level of reduced schemes. We do not know, whether the resulting
schemes actually coincide when we take the nilpotents into account.
And neither do we know whether the resulting reduced scheme is normal.

We should say that the problem of constructing the Uhlenbeck spaces
can be posed over a base field of any characteristic. However, the
proof  of one the main results of \cite{BrGa}, which insures that our spaces
$\calU^d_G$ are invariantly defined, uses the char.=0 assumption.
It is quite possible that in order to treat the char.=p case,
one needs a finer analysis.

\subsection{A sketch of the construction}
The construction of $\calU^d_G$ used in \cite{bfg} is a simplification
of a suggestion of Drinfeld's (the latter potentially  works for
an arbitrary surface $\bS$).
We are trying to express points of $\calU^d_G$ (one may call them quasi-bundles)
by replacing the original problem for the surface $\BP^2$ by another
problem for the curve $\BP^1$. Let us first generalize the problem to the case of
$G$-bundles with a parabolic structure along a fixed straight line.

Namely, let $\bfS=\PP^2$ and let $\PP^1_\infty\subset \bfS$ be the "infinite line"
(so that $\bfS\backslash \PP^1_\infty=\CC^2$). Let also
Let
$C\simeq \PP^1\subset\bfS$ denote the horizontal line in $\bfS$. Choose a parabolic
subgroup $P\subset G$. Let $\Bun_{G,P}$ denote the moduli space of
the following objects:

1) A principal $G$-bundle $\calH_G$ on $\bfS$;

2) A trivialization of $\calH_G$ on $\PP^1_\infty\subset \bfS$;

3) A reduction of $\calH_G$ to $P$ on $C$ compatible with the
trivialization of $\calH_G$ on $C$.

Let us describe the connected components of $\Bun_{G,P}$. Let $M$
be the Levi group of $P$. Denote by $\chM$ the {\it Langlands
dual} group of $M$ and let $Z(\chM)$ be its center. We denote by
$\Lam_{G,P}$ the lattice of characters of $Z(\chM)$. Let also
$\Lam_{G,P}^{\aff}=\Lam_{G,P}\x \ZZ$ be the lattice of characters
of $Z(\chM)\x\CC^*$. Note that $\Lam_{G,G}^{\aff}=\ZZ$.

The lattice $\Lam^{\aff}_{G,P}$ contains canonical semi-group
$\Lam^{\aff,+}_{G,P}$ of positive elements.
It is not difficult to see that the connected
components of $\Bun_{G,P}$ are parameterized  by the elements of
$\Lam_{G,P}^{\aff,+}$:
$$
    \Bun_{G,P}=\bigcup\limits_{\theta_{\aff}\in\Lam_{G,P}^{\aff,+}}
    \Bun_{G,P}^{\theta_{\aff}}.
$$

Typically, for $\theta_{\aff}\in \Lam_{G,P}^{\aff}$ we shall write
$\theta_\aff=(\theta,d)$ where $\theta\in \Lam_{G,P}$ and $d\in \ZZ$.

One would also like to construct the corresponding "Uhlenbeck scheme" $\calU_{G,P}^{\theta_{\aff}}$
stratified in the following way (the reader should compare it with \refe{st-fl}):
\eq{st-fl-aff}
\calU_{G,P}^{\theta}=\bigcup\limits_{\mu_{\aff}\in\Lam_{G,P}^{\aff,+}}\Bun_{G,P}^{\theta_{\aff}-\mu_{\aff}}\x
\Sym^{\mu_{\aff}}(\CC^2).
\end{equation}

The idea of the construction is as follows.
Let us consider  the scheme classifying triples
$(\F_G,\beta,\gamma)$, where

1) $\F_G$ is a principal $G$-bundle on $\PP^1$;

2) $\beta$ is a trivialization of $\F_G$ on the formal
neighborhood of $\infty\in\PP^1$;

3)  $\gamma$ is  a reduction to $P$ of the fiber of $\F_G$ at
$0\in\PP^1$.

It is easy to see that this scheme is canonically isomorphic to the thick
partial flag variety $\bfX_{G,P}^{\aff}=G(\calK)/\bfI_P$. Under this identification
the point $e_{G,P}^{\aff}\in\bfX_{G,P}^{\aff}$ corresponding to the unit element of $G$
corresponds to the trivial $\F_G$ with the trivial trivialization.

It is explained in \cite{bfg} that the variety
$\Bun_{G,P}$ is canonically isomorphic to the scheme classifying
{\it based maps} from $(C,\infty_C)$ to
$(\bfX_{G,P}^{\aff},e_{G,P}^{\aff})$ (i.e. maps from $C$
to $\bfX_{G,P}^{\aff}$ sending $\infty_C$ to
$e_{G,P}^{\aff}$).

One of the main results of \cite{bfg} gives an
explicit description of the Intersection Cohomology sheaf of all $\calU_{G,P}^{\theta_{\aff}}$.
We shall not reproduce the full answer here;
we shall only say that this answer is formulated in terms of the Lie algebra $\chg_{\aff}$ - the affine
Lie algebra whose Dynkin diagram is dual to that of $\grg_{\aff}$. Note that
in general $\chg_{\aff}\neq (\chg)_{\aff}$; in fact $\chg_{\aff}$ may result to be a twisted
affine Lie algebra (thus it is not isomorphic to the affinization of any finite-dimensional $\grg$).
We regard it as one of the first glimpses to (not yet formulated) Langlands duality for affine Lie algebras.

The proof of our computation of the IC-sheaves is also of independent interest. Namely, since in the affine
case the results of \cite{MV} are not available we must have a different way to see the algebra $\chg_{\aff}$
from the above geometry. In \cite{bfg} we first do on the combinatorial level; namely we realize the canonical
{\em Kashiwara crystal} discussed in \cite{KS} in terms of the varieties $\Bun_{G,B}^{\theta_{\aff}}$
(the idea of this realization is based on the earlier work \cite{BGcrys}). We then use this geometric construction
of crystals to compute the IC-sheaves of $\calU_{G,B}^{\theta_{\aff}}$ (the answer is very similar to Theorem \ref{IC})
which subsequently allows us to do it also for all $\calU_{G,P}^{\theta_{\aff}}$ using techniques
similar to those developed in \cite{BFGM} (in particular, we compute the IC-sheaf for $P=G$ which is
probably the most interesting case).

\subsection{The universal Verma module for $\chg_{\aff}$}\label{uhl-verma}
The scheme $\calU_{G,B}^{\theta_{\aff}}$ is endowed with a natural
action of $T\x (\CC^*)^2$ (here $T\subset G$ acts by changing the trivialization
of $\calH_G$ (cf. the the previous subsection) at $\PP^1_{\infty}$ and $(\CC^*)^2$
acts on $\bfS=\PP^2$ preserving $\PP^1_{\infty}$ and $C$).
Note that the field $\calK_{T\x (\CC^*)^2}$ can be thought of as a field of rational
functions of the variables $a\in \grt, \eps_1,\eps_2\in \CC$.
Define
$$
\IH_{G,B}^{\theta_{\aff}}=\IH^*_{T\x(\CC^*)^2}
(\calU_{G,B}^{\theta_{\aff}})\underset{\calA_{T\x(\CC^*)^2}}
\ten \calK_{T\x(\CC^*)^2}, \quad
\IH_{G,B}^{\aff}
=\bigoplus\limits_{\theta\in \Lam_{G,B}^{\aff,+}}\IH_{G,B}^{\theta_{\aff}}.
$$
Thus $\IH_{G,B}^{\theta_{\aff}}$ is a vector space over $\calK_{T\x(\CC^*)^2}$ which is endowed
with an intersection pairing $\la\cdot,\cdot\ra_{\theta_{\aff}}$ taking values in $\calK_{T\x(\CC^*)^2}$.
Also, for each $\theta_{\aff}$ as above we have the canonical element
$1_{G,B}^{\theta_{\aff}}\in\IH_{G,B}^{\theta_{\aff}}$ corresponding to the unit cohomology class.

In \cite{Br} we show that the Lie algebra $\chg_{\aff}$ acts naturally on $\IH_{G,B}^{\theta_{\aff}}$;
the corresponding $\chg_{\aff}$-module is naturally isomorphic to the Verma module
$M(\lam_{\aff})$ where $\lam_{\aff}=\frac{(a,\eps_2)}{\eps_1}+\rho_{\aff}$ (cf. \cite{Br} for more
details). Note that this is very similar to statement 1) from Subsection \ref{g-verma}.
We also have analogs of the statements 2), 3) 4) from Subsection \ref{g-verma}.

\section{Applications to gauge theory and quantum cohomology of
(affine) flag manifolds}
\label{nek}

\subsection{The partition function}
Recall that in the previous section we considered the moduli space $\Bun_G^d$ of $G$-bundles
on $\bfS=\PP^2$ trivialized at $\PP^1_{\infty}\subset \bfS$ and having 2nd Chern class equal to
$d$. We also have the scheme $\calU^d_G$ containing $\Bun_G^d$ as an open subset.
The group $G\x GL(2)$ acts naturally on $\calU^d_G$ where $G$ acts by changing the trivialization
at $\PP^1_{\infty}$ and $GL(2)$ acts on $\bfS$ preserving $\PP^1_{\infty}$.

 Thus we may consider (cf. \cite{Br}
\cite{nek}, \cite{nayo} for
precise definitions) the {\it equivariant integral}
$$
\int\limits_{\calU^d_G}1^d
$$
of the unit $G\times GL(2)$-equivariant cohomology class (which we
denote by $1^d$) over
$\calU_G^d$; the integral takes values in the field $\calK$ which is
the field of fractions of the algebra $\calA=H^*_{G\x
GL(2)}(pt)$. Note that $\calA$ is
canonically isomorphic to the algebra of polynomial functions on the Lie algebra
$\grg\x\gl(2)$  which
are invariant with respect to the adjoint action.
Thus each $\int\limits_{\calU^d_G}1^d$ may naturally be regarded as a
rational function of $a\in\grt$ and $(\eps_1,\eps_2)\in \CC^2$; this function
must be invariant with respect to the natural action of $W$ on $\grt$ and with respect
to interchanging $\eps_1$ and $\eps_2$.

Consider now the generating function
$$
\calZ=\sum\limits_{d=0}^\infty Q^d \int\limits_{\calU_G^d} 1^d.
$$
It can (and should) be thought of as a function of the variables
$\grq$ and $a,\eps_1,\eps_2$ as before. In \cite{nek} it was
conjectured that the first term of the asymptotic in the limit
$\underset{\eps_1,\eps_2\to 0}\lim \ln \calZ$ is closely related to
{\it Seiberg-Witten prepotential} of $G$.
\footnote{In fact, in \cite{nek} this conjecture is only formulated for $G=SL(n)$ but
the generalization to other groups is pretty straightforward. Also, one can reformulate
everything in the language of moduli spaces of anti-selfdual connections of the 4-sphere
$S^4$ rather than in terms of holomorphic (=algebraic) G-bundles on $\PP^2$; this, perhaps,
is closer to the physical origins of the problem.}
This can be thought of as a rigorous
mathematical formulation of the results of Seiberg and Witten from 1994.
For $G=SL(n)$ this
conjecture has been proved in \cite{neok} and \cite{nayo}. Also in
\cite{nek} an explicit combinatorial expression for $\calZ$ has
been found. We are going to give a sketch of the proof of this conjecture for arbitrary $G$.
%---------------------------------------------------------------------------------------------------------------
\subsection{Parabolic generalization of the partition function}
Recall from the previous section that for any parabolic $P\subset G$ we have the varieties
$\Bun_{G,P}^{\theta_{\aff}}$ and $\calU_{G,P}^{\theta}$. The latter contains the former
as a dense open subset. Is is easy to see that $\calU_{G,P}^{\theta}$ has a natural
action of the group $M\x (\CC^*)^2$ where $M$ denotes the Levi subgroup of $P$.
Also we have the field $\calK_{M\x (\CC^*)^2}$ which is isomorphic to the field of rational functions on
$\grm\x\CC^2$ which are invariant with respect to the adjoint action.

Let $T\subset M$ be a maximal torus. Then one can show that
$(\calU_{G,P}^{\theta_{\aff}})^{T\x(\CC^*)^2}$ consists of one point.
This guarantees that we may consider the integral
$\int\limits_{\calU_{G,P}^{\theta_{\aff}}}1_{G,P}^{\theta_{\aff}}$ where
$1_{G,P}^{\theta_{\aff}}$
denotes the unit class in
$H^*_{M\x(\CC^*)^2}(\calU_{G,P}^{\theta_{\aff}},\CC)$.
The result can be thought
of as a rational function on $\grm\x\CC^2$ which is invariant with
respect to the adjoint action of $M$. Define
\eq{partition-affine}
    \calZ_{G,P}^{\aff}=\sum\limits_{\theta\in\Lam_{G,P}^{\aff}} \grq_{\aff}^{\theta_{\aff}}
    \ \int\limits_{\calU_{G,P}^{\theta_{\aff}}}1_{G,P}^{\theta_{\aff}}.
\end{equation}
One should think of $\calZ_{G,P}^{\aff}$ as a formal power series
in $\grq_{\aff}\in Z(\chM)\x\CC^*$ with values in the space of ad-invariant
rational functions on $\grm\x\CC^2$. Typically, we shall write
$\grq_{\aff}=(\grq,Q)$ where $\grq\in Z(\chM)$ and $Q\in\CC^*$.
Also we shall denote an element of $\grm\x\CC^2$ by
$(a,\eps_1,\eps_2)$ or (sometimes it will be more convenient) by
$(a,\hbar,\eps)$ (note that for general $P$ (unlike in the case $P=G$)
the function $\calZ_{G,P}^{\aff}$ is not symmetric with respect to
switching
$\eps_1$ and $\eps_2$).

%-----------------------------------------------------------
\subsection{The "finite-dimensional" analog}
Recall that the space $\calU_{G,P}^{\theta_{\aff}}$ is closely related to the space
of based quasi-maps $C=\PP^1\to \bfX_{G,P}^{\aff}$ of degree $\theta_{\aff}$.
Since the scheme $\bfX_{G,P}^{\aff}$ may (and should) be thought of as
a partial flag variety for $\grg_{\aff}$ it is natural to consider the
following "finite-dimensional" analog of the above problem.
Recall that we denote by $^{\bt}\calM_{G,P}^{\theta}$  the moduli space of based
maps from $(C,\infty_{C})$ to
$(X_{G,P},e_{{G,P}})$ of degree $\theta$, i.e. the moduli space of maps
$C\to X_{G,P}$ which send $\infty_{C}$ to $e_{G,P}$.
This space is acted on by the group $M\x\CC^*$.

We now introduce the ``finite-dimensional'' analog of the partition
function
\refe{partition-affine}. As before let
$$
\calA_{M\x\CC^*}=H^*_{M\x\CC^*}(pt,\CC)
$$
and  denote by $\calK_{M\x\CC^*}$ its field of fractions. Let also
$1_{G,P}^{\theta}$ denote the
unit class in the
$M\x \CC^*$-equivariant cohomology of $^{\bt}\QM_{G,P}^{\theta}$.
Then we define
\eq{partition-finite}
       \calZ_{G,P}=\sum\limits_{\theta\in\Lam_{G,P}^{\theta}}
\grq^{\theta} \ \int\limits_{^{\bt}\QM_{G,P}^{\theta}}1_{G,P}^{\theta}.
\end{equation}
This is a formal series in $q\in Z(\chM)$ with values in the field
$\calK_{M\x\CC^*}$ of $M$-invariant rational functions on $\grm\x\CC$.

In fact the function $\calZ_{G,P}$ is a familiar object in Gromov-Witten theory: in
\cite{Br} we show
 that up to a simple factor $\calZ_{G,P}$
is the so called {\it equivariant $J$-function} of $X_{G,P}$; this function is in some
sense responsible for the (small) quantum cohomology of $X_{G,P}$. Thus the problem of
computation of the function
$\calZ_{G,P}^{\aff}$ may be thought of as the problem of computation of the (not yet rigorously defined)
quantum cohomology of the affine flag manifolds $\bfX_{G,P}^{\aff}$. Note that in the case
$G=SL(n)$ and $P=B$ a heuristic computation of the latter ring in terms of the so called
{\em periodic Toda lattice} was done in \cite{go}; our results discussed presented in the nex
subsection are compatible with this computation.
%---------------------------------------------------------
\ssec{borel-int}{Computation of the partition functions in the Borel case}
We believe that it should be possible to
express the function $\calZ_{G,P}$ (resp. the function
$\calZ_{G,P}^{\aff}$) in terms of representation theory of the Lie
algebra $\chg$ (resp. $\chg_{\aff}$) -- by the definition this is a
Lie algebra whose root system is dual to that of $\grg$ (resp. to that
of $\grg_{\aff}$).  One of the main
results of \cite{Br} gives such a calculation of the functions
$\calZ_{G,B}$ and $\calZ_{G,B}^{\aff}$ where $B\subset G$ is a
Borel subgroup of $G$. Roughly speaking we show that $\calZ_{G,B}$
(resp. $\calZ_{G,B}^{\aff}$) is equal to {\it Whittaker matrix
coefficient} of the Verma module over $\chg$ (resp. over
$\chg_{\aff}$) whose lowest weight given by $\frac{a}{\hbar}+\rho$
(resp. $\frac{(a,\eps_1)}{\eps_2}+\rho_{\aff}$ where $a,\hbar,\eps_1$ and $\eps_2$
are as in Subsection \label{uhl-verma} (here we regard $(a,\eps_1)$ as a
weight
for the dual affine algebra $\chg_{\aff}$; this is explained carefully
in Section 3 of \cite{Br}). These statements in fact follow immediately
from the results of subsections \ref{g-verma} and \ref{uhl-verma}
after one formulates the definition of equivariant integration using
intersection cohomology (this is done is \cite{Br}).

The above description of the partition function allows one to
produce certain differential equations which are satisfied by
the functions $\calZ_{G,B}$ and $\calZ_{G,B}^{\aff}$. More precisely,
we show that the function
$\grq^{\frac{a}{\hbar}}\calZ_{G,B}$ is an eigen-function of
the {\em quantum Toda hamiltonians}
 associated with $\chg$
with eigen-values determined (in the natural way) by $a$
 (we refer the reader to
\cite{etingof} for the definition of (affine) Toda integrable system
and its relation with Whittaker functions). In this way we reprove the
results of \cite{gk} and \cite{kim} about (equivariant) quantum cohomology
of the flag varieties $X_{G,P}$.
In the affine case one can also show that
$\grq^{\frac{a}{\hbar}}\calZ_{G,B}^{\aff}$
is an eigen-function of a certain differential operator which has
order 2 (``non-stationary analog'' of the affine quadratic Toda hamiltonian).
In \cite{BrEt} we explain how this allows to compute the asymptotics of {\it all}
the functions $\calZ_{G,P}^{\aff}$ when $\eps_1,\eps_2\to 0$ in
another publication. We also show in \cite{BrEt} that this implies the Nekrasov conjecture
mentioned above for arbitrary $G$.

\section{Some open problems}\label{open}

In this section we present a list of open problems that  related to the subjects in the
preceding sections of the paper.

\subsection{Normality and rational singularities}
We conjecture that the schemes $\QM_{G,P}^{\theta}$ and $\calU_{G,P}^{\theta_{\aff}}$ are normal and have rational
singularities. From purely technical point of view one needs to know this in order to attack problem 2. On the other hand,
this statement seems to be important for the following reason. It is explained in [9] and [14] that one should think
about the singularities of the schemes $\QM_{G,P}^{\theta}$ as a finite-dimensional model for the singularities
of the so called {\it semi-infinite Schubert varieties} (cf. [14] for a more detailed discussion of this).
On the other hand, one knows (cf. [13] and references therein) that usual (both finite and affine) Schubert varieties
{\it are} normal and have rational singularities. Hence in the case of $\QM_{G,P}^{\theta}$ our conjectures
can be thought of as a generalization of this result to the semi-infinite Schubert varieties.

\subsection{Computation of parabolic partition functions}It would be very interesting to
express the functions $\calZ_{G,P}$ (resp. $\calZ_{G,P}^{\aff}$) in terms of representation theory of the algebra
$\chg$ (resp. $\chg_{\aff}$). Let us note that a combinatorial expression for all these functions in the case
$G=SL(n)$ (in the finite case) was found in \cite{yau}, but we do not know how to interpret this answer in terms
of representation theory.

\subsection{Other cohomology theories}
The functions $\calZ_{G,P}$ and $\calZ_{G,P}^{\aff}$ have their $K$-theoretic counterparts
$\calZ^K_{G,P}$ and $\calZ^{K,\aff}_{G,P}$ (to define those one needs
to replace the equivariant integrals considered in \refe{partition-affine} and \refe{partition-finite} by
the corresponding integrals in equivariant $K$-theory). The function $\calZ^K_{G,P}$ is exactly the $K$-theoretic $J$-function
of $X_{G,P}$ as defined in [22]. We would like to express these functions using the representation
theory of the quantum group $U_q(\chg)$ (resp. $U_q(\chg_{\aff})$). Presumably, in order to do this one should
be able to construct geometrically some representations of these quantum groups.
For $P=B$ and $G=SL(n)$ this is done in [7]. For $P=G=SL(n)$ the $K$-theoretic partition function $\calZ_G^{K,\aff}$
was also studied in \cite{neok}.

Let us also recall that in [6] the authors use the results described in \refss{borel-int} in order to connect certain
asymptotic of the function $\calZ_{G,P}^{\aff}$ with the {\it Seiberg-Witten prepotential} corresponding to
the classical Toda integrable system associated with the Lie algebra $\chg_{\aff}$. We would like to generalize this
to the $K$-theoretic partition functions. This should involve some interesting interplay between quantum affine
algebras and certain difference equations (which should be thought of non-stationary deformations of known
integrable difference equations of Toda type). The corresponding classical integrable system describing the
asymptotic should be the so called {\it relativistic Toda system} associated with the Lie algebra $\chg_{\aff}$.

It would also be interesting to generalize this to the case when $K$-theory is replaced by any elliptic cohomology
theory. Again, for $P=G=SL(n)$ (in the affine case) this is done in \cite{vafa}.

\subsection{Chern classes of the tangent bundle and the Calogero-Moser system}
As was mentioned above the partitions functions $\calZ_{G,P}^{\aff}$ are related to
the pure N=2 super-symmetric gauge theory in 4 dimensions. One should also have an extension
of the above result for gauge theory {\it with matter}. This means that instead of considering equivariant
integrals of the unit cohomology class we should consider integrals of the Chern classes of various natural
bundles on the moduli spaces in question. For example in the case of {\it adjoint matter} one should
integrate the Chern polynomial of the tangent bundle of $\calU_{G,P}^{\theta_{\aff}}$ (this, of course,
has to be properly interpreted since the variety in question is singular). For $P=G=SL(n)$ such functions
are studied in [28] and the corresponding asymptotic of the the partition function is shown there
to be related with the prepotential of the classical {\it elliptic Calogero-Moser} integrable system.
We expect that for $P=B$ (and for general $G$) the corresponding partition function should be closely
related with the universal eigen-function of the corresponding non-stationary deformation of the quantum
Calogero-Moser hamiltonian associated with the Lie algebra $\chg_{\aff}$. Similar statement should also hold
in the finite case (i.e. when we integrate over $\QM_{G,B}^{\theta}$'s and not over $\calU_{G,B}^{\theta_{\aff}}$'s).
In particular, in the finite case we should get a geometric interpretation of the universal eigen-function of the
quantum Calogero-Moser system associated with the Lie algebra $\chg$.

\subsection{Functional equation for parabolic Eisenstein series}
It will be very important to generalize the functional equation for geometric Eisenstein series
discussed in Subsection \ref{func} to the case of parabolic Eisenstein series. More precisely,
given a parabolic subgroup $P\subset G$ with the Levi subgroup $M$ in Section \ref{eisen} we discussed
the Eisenstein series functor $\Eis_M^G$. In fact this notation is slightly misleading since this
functor actually depends not just on $M$ but also on $P$. Let us now use the notation $\Eis_{M,P}^G$
for it. With this notation the "functional equation" problem reads as follows:
given two parabolic subgroups $P,Q\subset G$ containing {\em the same} Levi subgroup $M$, construct
as isomorphism between the functors $\Eis_{M,P}^G$ and $\Eis_{M,Q}^G$ restricted to some large
subcategory of "regular" sheaves inside $\Sh(M)$.

\end{document}